\documentclass[10pt]{article}%
\usepackage{latexsym,amsmath,amssymb,graphics,amscd}
\usepackage{amsmath}
\usepackage{amsfonts}
\usepackage{amssymb}
\usepackage{graphicx}%
\newcommand{\bfi}{\bfseries\itshape}
\addtocounter{footnote}{1}
\makeatletter
\@addtoreset{figure}{section}
\def\thefigure{\thesection.\@arabic\c@figure}
\def\fps@figure{h,t}
\@addtoreset{table}{bsection}
\def\thetable{\thesection.\@arabic\c@table}
\def\fps@table{h, t}
\@addtoreset{equation}{section}

\makeatother
\newtheorem{theorem}{Theorem}[section]
\newtheorem{definition}[theorem]{Definition}
\newtheorem{lemma}[theorem]{Lemma}
\newtheorem{remark}[theorem]{Remark}
\newtheorem{proposition}[theorem]{Proposition}
\newtheorem{corollary}[theorem]{Corollary}

\newsavebox{\savepar}

\makeatletter
\begin{document}

\title{\textbf{Superposition rules and stochastic Lie-Scheffers systems}}
\author{Joan-Andreu L\'{a}zaro-Cam\'{\i}$^{1}$ and Juan-Pablo Ortega$^{2}$}
\date{}
\maketitle

\begin{abstract}
This paper proves a version for stochastic differential equations of the Lie-Scheffers Theorem. This result characterizes the existence of nonlinear superposition rules for the general solution of those equations in terms of the involution properties of the distribution generated by the vector fields that define it. When stated in the particular case of standard deterministic systems, our main theorem improves various aspects of the classical Lie-Scheffers result. We show that the stochastic analog of the classical Lie-Scheffers systems can be reduced to the study of Lie group valued stochastic Lie-Scheffers systems; those systems, as well as those taking values in homogeneous spaces are studied in detail. The developments of the paper are illustrated with several examples. 
\end{abstract}

\makeatletter
\addtocounter{footnote}{1} \footnotetext{Departamento de F\'{\i}sica
Te\'{o}rica. Universidad de Zaragoza. Pedro Cerbuna, 12. E-50009 Zaragoza.
Spain. {\texttt{lazaro@unizar.es}}} \addtocounter{footnote}{1}
\footnotetext{Centre National de la Recherche Scientifique (CNRS), D\'{e}partement de
Math\'{e}matiques de Besan\c{c}on, Universit\'{e} de Franche-Comt\'{e}, UFR
des Sciences et Techniques. 16, route de Gray. F-25030 Besan\c{c}on cedex.
France. {\texttt{Juan-Pablo.Ortega@univ-fcomte.fr} }} \makeatother

\section{Introduction}

A differential equation is said to have a {\bfi  superposition rule} (a more explicit definition is provided in the next section) whenever any of its solutions can be written as a given (in general nonlinear) function of the initial condition and of a fixed set of particular solutions. The first characterization of the existence of superposition rules was given by the Norwegian mathematician Sophus Lie in a remarkable piece of work~\cite{lie 1893} where he established a link between the existence of superposition rules and what we nowadays call the Lie algebraic properties of the vector fields that define a time-dependent differential equation. This result is referred to as the {\bfi  Lie-Scheffers Theorem} and systems that satisfy its hypotheses as {\bfi  Lie-Scheffers systems}.

Lie-Scheffers systems have been the subject of much attention due to their widespread occurrence in physics and mathematics. The reader is encouraged to check with~\cite{libro marmo, Cari07}, and references therein, for various presentations of the classical Lie-Scheffers Theorem, an excellent collection of examples of applications of this theorem, and for historical remarks.

The main goal of this paper is the extension of the Lie-Scheffers Theorem to stochastic differential equations. This generalization is stated in Theorem \ref{teorema Lie}. It is worth emphasizing that the main result of the paper,~Theorem \ref{teorema Lie}, cannot be seen just as a mere transcription of the
deterministic Lie-Scheffers Theorem into the context of Stratonovich stochastic integration by using the so called  Malliavin's Transfer
Principle~\cite{malliavin transfer}. This Principle states that whatever is true for standard differential equations also holds for Stratonovich stochastic differential equations; as we will see later on, there are purely stochastic conditions that appear in the statement of the theorem. 

Additionally, in proving Theorem \ref{teorema Lie}  we have  carefully spelled out the regularity conditions needed for the result to be valid; those conditions are only vaguely evoked in the classical references or in the cited papers that  study the deterministic case. More importantly, a careful construction of the proof has lead us to realize that the hypotheses under which we can guarantee the existence of superposition rules can be weakened: the Lie algebra condition in the classical theorem can be replaced by an involutivity hypothesis that is, in general, less restrictive. 

The contents of the paper are structured as follows. Section~\ref{Superposition rules for stochastic differential equations} explains in detail the notion of superposition rule and includes a proposition that translates this concept into geometric terms. Section~\ref{The stochastic Lie-Scheffers Theorem} contains the main theorem that we have already described. 

Section~\ref{Lie-Scheffers systems on Lie groups and homogeneous spaces} is dedicated to the study of Lie-Scheffers systems on Lie groups and homogeneous spaces; this case is particularly relevant since, as we show in the first result of that section (Proposition~\ref{locally on lie groups}), classical Lie-Scheffers systems (roughly speaking, those generated by vector fields that close a Lie algebra) can be locally reduced to this case via a theorem due to Palais. In that section we also show, as an example, how L\'evy stochastic processes can be seen as Lie group valued Lie-Scheffers systems. The section concludes with a brief presentation of the classical Wei-Norman method for solving Lie-Scheffers systems, adapted to the stochastic context.

Section~\ref{The flow of a stochastic Lie-Scheffers system} contains a discussion on how the existence of a superposition rule for a stochastic differential equation makes available a remarkable feature that has deserved certain  attention in the context of standard stochastic differential equations, namely, the fact that the stochastic flow can be written as a fixed deterministic function of the Brownian forcing of the equation in question. Indeed, a well know theorem by Ben Arous~\cite{Ben Arous}, that we state in the paper and whose proof is based on the use of stochastic Taylor expansions, shows that this property of the flow is available under exactly the same hypotheses as the classical Lie-Scheffers Theorem. Our main theorem allows, admittedly only to a certain extent, the generalization of this statement to any stochastic differential equation that satisfies its hypotheses; more specifically, any SDE generated by vector fields that span an involutive distribution has a superposition rule and hence its flow can be written as a fixed deterministic function of the initial conditions and of a set of solutions that contain the stochastic behavior of the resulting map.

The paper concludes with a section that contains a number of examples that illustrate the developments of the paper.

\section{Superposition rules for stochastic differential equations}
\label{Superposition rules for stochastic differential equations}

Let $(\Omega, \mathcal{F}, P)$ be a probability space. We start by considering the stochastic differential equation
\begin{equation}
\label{stochastic differential equation expression}\delta\Gamma=S\left(
X,\Gamma\right)  \delta X,
\end{equation}
where $X:\mathbb{R}_{+}\times\Omega\rightarrow\mathbb{R}^{l}$ is a given
$\mathbb{R}^{l}$-valued semimartingale and $S\left(  x,z\right)
:T_{x}\mathbb{R}^{l}\longrightarrow T_{p}\mathbb{R}^{n} $ is a Stratonovich
operator from $\mathbb{R}^{l}$ to $\mathbb{R}^{n}$. Sometimes we will choose a basis in $T ^\ast \mathbb{R} ^l $ and will write down the Stratonovich operator $S(x,z)$ in terms of its components $(S_{1}\left(x,z\right)  ,\ldots,S_{l}(x,z))$ with respect to that basis.

\begin{definition}
\label{def 1} A {\bfi superposition rule} of the stochastic differential
equation~(\ref{stochastic differential equation expression}) is a pair
$(\Phi,\{\Gamma_{1},\ldots,\allowbreak\Gamma_{m}\})$, where $\Phi
:\mathbb{R}^{n\left(  m+1\right)  }\longrightarrow\mathbb{R}^{n}$ is a (not
necessarily smooth) function and $\{\Gamma_{i}:\mathbb{R}_{+}\times
\Omega\rightarrow\mathbb{R}^{n}~|~i=1,\ldots,m\}$ is a set of particular
solutions of~(\ref{stochastic differential equation expression}) such that any
solution $\Gamma$ of~(\ref{stochastic differential equation expression}) can
be written, at least up to a sufficiently small stopping time $\tau$, as
\[
\Gamma=\Phi\left(  z^{1},\ldots,z^{n};\Gamma_{1},\ldots,\Gamma_{m}\right)
=:\Phi\left(  z;\Gamma_{1},\ldots,\Gamma_{m}\right)  ,
\]
where $z=\left(  z^{1},\ldots,z^{n}\right)  $ a set of $n$ arbitrary constants
associated with the initial condition of the solution $\Gamma$, that is,
$\Gamma(0,\omega)=(z^{1},\ldots,z^{n})$, for all $\omega\in\Omega$. We extend
to the stochastic context the terminology used for standard differential
equations and we will call {\bfi Lie-Scheffers} systems the stochastic
differential equations that admit a superposition rule.
\end{definition}

\begin{remark}
\label{only euclidean} \normalfont As we will see in examples later on in the
paper, superposition rules exist only locally. That is why we can, without
loss of generality, restrict our attention to stochastic differential equation
on Euclidean spaces. Observe also  that we are requiring that $\Phi$ does not depend on
time, the probability space, or the noise $X$. This prevents us from using certain
regularization techniques at the time of testing the existence of superposition rules.
For example, when dealing with a deterministic differential equation, the standard
transformation of a time-dependent system $\dot{\gamma}=f\left(
t,\gamma\right)  $ on $\mathbb{R}^{n}$, $f:\mathbb{R}^{n+1}\rightarrow
\mathbb{R}^{n}$ into the autonomous one
\[
\dot{\gamma}=f\left(  t,\gamma\right)  \text{ \ \ and \ \ }\dot{t}=1
\]
on $\mathbb{R}^{n+1}$ obtained by adding an extra trivial differential equation for the
time, is not allowed; indeed, if we find a superposition rule for the transformed autonomous system, that rule does not yield a superposition rule for the original system that satisfies the requirements of our definition, precisely due to the explicit dependence on time that appears in the superposition  function.
\end{remark}

In order to study the implications of the presence of a superposition rules we
take a more geometric approach. Let $\Psi$ be the function defined by
\begin{equation}%
\begin{array}
[c]{rrl}%
\Psi:\mathbb{R}^{n\left(  m+2\right)  } & \longrightarrow & \mathbb{R}^{n}\\
\left(  z,q_{0},q_{1},\ldots,q_{m}\right)  & \longmapsto & q_{0}-\Phi\left(
z;q_{1},\ldots,q_{m}\right)  .
\end{array}
\label{eq 2}%
\end{equation}
Notice that for any $z\in\mathbb{R}^{n}$, the function $\Psi_{z}:=\Psi\left(
z,\cdot\right)  :\mathbb{R}^{n\left(  m+1\right)  }\rightarrow\mathbb{R}^{n}$
is constant on a $\left(  m+1\right)  $-tuple $\left(  \Gamma,\Gamma_{1}%
\ldots,\Gamma_{m}\right)  $ of solutions of the system
(\ref{stochastic differential equation expression}), at least up to a given
stopping time $\tau$, provided that $\Gamma_{t=0}=z\in\mathbb{R}^{n}$ a.s..
From now on we assume that all the solutions $\Gamma$ that we are dealing with are
constant a.s. at $t=0$. Additionally, if the function $\Phi$ is smooth then
the map $\Psi_{z}:\mathbb{R}^{n\left(  m+1\right)  }\rightarrow\mathbb{R}^{n}$
is a submersion for any fixed $z\in\mathbb{R}^{n}$, because
\begin{equation}
\operatorname{rank}\left(  \frac{\partial\Psi_{z}^{j}}{\partial q_{0}^{i}%
}\right)  _{j,i=1,\ldots,n}=\operatorname{rank}\left(  I_{n}\right)  =n
\label{eq 18}%
\end{equation}
where $I_{n}$ is the identity matrix of dimension $n$. Consequently, for any
$z\in\mathbb{R}^{n}$, the level set $\Psi_{z}^{-1}\left(  0\right)
\subset\mathbb{R}^{n\left(  m+1\right)  }$ is a closed embedded submanifold of
$\mathbb{R}^{n\left(  m+1\right)  }$ of dimension $nm$ . That is, the function
$\Psi$ defines a family $\mathcal{G}$ of regular $nm$-dimensional submanifolds
$\mathcal{G}_{z}$ via the zero level sets $\Psi_{z}^{-1}\left(  0\right)
=\{p\in\mathbb{R}^{n(m+1)}\mid\Psi\left(  z,p\right)  =0\}=:\mathcal{G}_{z}$
of $\Psi_{z}$, for any $z\in\mathbb{R}^{n}$. The submanifolds $\mathcal{G}%
_{z}$ are globally diffeomorphic to $\mathbb{R}^{nm}$ via the restriction
$\pi_{m}|_{\mathcal{G}_{z}}$ to $\mathcal{G}_{z}$ of the projection $\pi
_{m}:\mathbb{R}^{n(m+1)}=\mathbb{R}^{n}\times\overset{m+1}{\cdots}%
\times\mathbb{R}^{n}\longrightarrow\mathbb{R}^{nm}=\mathbb{R}^{n}%
\times\overset{m}{\cdots}\times\mathbb{R}^{n}$ onto the last $m$
$\mathbb{R}^{n}$ factors. This is easy to see by verifying that the inverse
$\Xi_{z}:\mathbb{R}^{mn}\rightarrow\mathcal{G}_{z}$ of $\pi_{m}|_{\mathcal{G}%
_{z}}$ is given by $\Xi_{z}(q_{1},\ldots,q_{m})=(\Phi(z;q_{1},\ldots
,q_{m}),q_{1},\ldots,q_{m})$, which is obviously a diffeomorphism. In order to
study the significance of the family of submanifolds $\mathcal{G}$ we start by
introducing the following definition.

\begin{definition}
\label{diagonal extension} Let $Y:\mathbb{R}^{n}\rightarrow\mathbb{R}^{n}$ be
a vector field. The vector field
\[%
\begin{array}
[c]{rrl}%
\widetilde{Y}:\mathbb{R}^{n\left(  m+1\right)  } & \longrightarrow &
\mathbb{R}^{n\left(  m+1\right)  }\\
\left(  q_{0},\ldots,q_{m}\right)  & \longmapsto & \left(  Y\left(
q_{0}\right)  ,\ldots,Y\left(  q_{m}\right)  \right)
\end{array}
\]
is called the \textbf{diagonal extension} of $Y$.
\end{definition}

It can be easily checked that the set of diagonal extensions of vector fields
in $\mathfrak{X}(\mathbb{R}^{n})$ are a subalgebra of $\mathfrak{X}%
(\mathbb{R}^{n(m+1)})$; more explicitly, for any $Y_{1}, Y_{2},Y_{3}%
\in\mathfrak{X}(\mathbb{R}^{n})$ and $\lambda\in\mathbb{R} $,
\begin{equation}
\label{eq 14}\lbrack\widetilde{Y}_{1},\widetilde{Y}_{2}+\lambda\widetilde
{Y}_{3}]=\widetilde{\left[  Y_{1} ,Y_{2}+\lambda Y _{3}\right]  }.
\end{equation}

\noindent The following proposition states that, roughly speaking, the family
of submanifolds $\mathcal{G}$ completely characterizes the superposition rule.

\begin{proposition}
\label{prop 2} Suppose that the stochastic differential equation
(\ref{stochastic differential equation expression}) admits a smooth
superposition rule $( \Phi,\left\{  \Gamma_{1},\ldots,\Gamma_{m}\right\}  ) $.
Suppose that $( \Gamma_{1},\ldots,\Gamma_{m}) _{t=0}=(p _{1}, \ldots, p _{m})
\in\mathbb{R}^{mn}$ a.s.. Then, there exists a family $\mathcal{G} $ of closed
embedded $nm $-dimensional submanifolds of $\mathbb{R}^{n\left(  m+1\right)
}$ such that for any $z \in\mathbb{R}^{n} $ there exists $\mathcal{G}_{z}
\in\mathcal{G} $ such that $( \Gamma^{z}, \Gamma_{1},\ldots,\Gamma_{m}%
)\subset\mathcal{G} _{z}$, with $\Gamma^{z} $ the solution
of~(\ref{stochastic differential equation expression}) such that $(\Gamma
^{z})_{t=0}=z$. Moreover, for any $\mathcal{G} _{z}\in\mathcal{G} $ the map
$\pi_{m}|_{\mathcal{G}_{z}}: \mathcal{G}_{z} \rightarrow\mathbb{R}^{nm}$ is a diffeomorphism.

Conversely, let $\mathcal{G}$ be a family of (not necessarily embedded)
submanifolds of $\mathbb{R}^{n\left(  m+1\right)  }$ diffeomorphic to
$\mathbb{R}^{nm}$ via $\pi_{m}$ and $\left\{  \Gamma_{1},\ldots,\Gamma
_{m}\right\}  $ a set of distinct solutions of
(\ref{stochastic differential equation expression}) such that $(\Gamma
_{1},\ldots,\Gamma_{m})_{t=0}=(p_{1},\ldots,p_{m})\in\mathbb{R}^{mn}$ a.s..
Then, if for any point $z\in\mathbb{R}^{n}$ there is an element $\mathcal{G}%
_{z}$ that contains the point $(z,p_{1},\ldots,p_{m})$ and the diagonal
extensions $(\widetilde{S}_{1}\left(  X,\cdot\right)  ,\ldots,\widetilde
{S}_{l}\left(  X,\cdot\right)  )$ of the vector fields $(S_{1}\left(
X,\cdot\right)  ,\ldots,S_{l}(X,\cdot))$ that
define~(\ref{stochastic differential equation expression}) are tangent to
$\mathcal{G}_{z}$ when evaluated at $(\Gamma^{z},\Gamma_{1},\ldots,\Gamma
_{m})$, then (\ref{stochastic differential equation expression}) admits a
(possibly nonsmooth) superposition rule.
\end{proposition}

\noindent\textbf{Proof.\ \ } In view of the remarks preceding
Definition~\ref{diagonal extension} we just need to prove that having a family
$\mathcal{G}$ that satisfies the hypotheses in the statement allows us to
recover the superposition rule.

Let $\left\{  \Gamma_{1},\ldots,\Gamma_{m}\right\}  $ be the set of fixed
distinct solutions of~(\ref{stochastic differential equation expression}).
Denote $p_{i}=(\Gamma_{i})_{t=0}$ the (necessarily different) constant initial
conditions of $\Gamma_{i}$, $i=1,\ldots,m$. Let $z=\left(  z^{1},\ldots
,z^{n}\right)  \in\mathbb{R}^{n}$ be a point and let $\mathcal{G}_{z}$ be the
submanifold in $\mathcal{G}$ such that $\left(  z,p_{1},\ldots,p_{m}\right)
\in\mathcal{G}_{z} $; by hypothesis, this manifold is diffeomorphic to
$\mathbb{R} ^{nm}$ via the map $\varphi_{z}=\left.  \pi_{m}\right\vert
_{\mathcal{G}_{z}}$, where $\pi_{m}:\mathbb{R}^{n\left(  m+1\right)
}\longrightarrow\mathbb{R}^{nm}$ is the projection onto the last $nm$ factors.
In other words, the last $nm$ coordinates of a point in $\mathbb{R}^{n\left(
m+1\right)  }$ serve as global coordinates of $\mathcal{G}_{z}$. Introduce the
projection%
\begin{equation}%
\begin{array}
[c]{rrl}%
\pi_{\mathbb{R}^{n}}^{0}:\mathbb{R}^{n(m+1)} & \longrightarrow &
\mathbb{R}^{n}\\
\left(  q_{0},\ldots,q_{m}\right)  & \longmapsto & q_{0}.
\end{array}
\label{eq 15}%
\end{equation}
We now define
\begin{equation}
\left(  \Gamma_{0}\right)  _{t}\left(  \omega\right)  :=\pi_{\mathbb{R}^{n}%
}^{0}\circ\varphi_{z}^{-1}\left(  \left(  \Gamma_{1}\right)  _{t}\left(
\omega\right)  ,\ldots,\left(  \Gamma_{m}\right)  _{t}\left(  \omega\right)
\right)  . \label{eq 11}%
\end{equation}
It is immediate to see that $\left(  \Gamma_{0}\right)  _{t=0}=z$ and that
$\Gamma_{0}$ is a semimartingale because, by construction, it is a composition
of smooth functions with semimartingales. Let now $\Gamma^{z}$ be the unique
solution of (\ref{stochastic differential equation expression}) with a.s.
initial condition $z\in\mathbb{R} ^{n}$. We will proceed by proving that
$\Gamma_{0}$ defined in (\ref{eq 11}) equals $\Gamma^{z}$ and we will
therefore have a superposition rule $\Phi$ given by the map $\Phi(z;
\Gamma_{1}, \ldots, \Gamma_{m}):=\pi_{\mathbb{R}^{n} }^{0}\circ\varphi
_{z}^{-1}\left(  \Gamma_{1}, \ldots, \Gamma_{m}\right)  $. Notice that unless
additional hypotheses are assumed on the family $\mathcal{G} $, there is no
guarantee on the smoothness of $\Phi$ on the $z$ variable.

In order to prove that $\Gamma_{0}$ equals $\Gamma^{z}$, denote by $\left(
q^{k};k=1,\ldots,n\right)  $ the coordinates on $\mathbb{R} ^{n}$ and by
$\left(  q_{a}^{k};k=1,\ldots,n;a=0,\ldots,m\right)  $ the coordinates on
$\mathbb{R}^{n(m+1)}$. Let $F_{k}^{a}:\mathbb{R}^{nm}\rightarrow\mathbb{R}%
^{n}$ and $X_{k}^{a}:\mathbb{R}^{nm}\rightarrow\mathbb{R}^{n(m+1)}$ be the
maps defined as%
\[
F_{k}^{a}\left(  q_{1},\ldots,q_{m}\right)  =T_{\left(  q_{1},\ldots
,q_{m}\right)  }(\pi_{\mathbb{R}^{n}}^{0}\circ\varphi_{z}^{-1}\circ\pi
_{m})\left(  \frac{\partial}{\partial q_{a}^{k}}\right)
\]%
\[
X_{k}^{a}\left(  \varphi^{-1}_{z}( q_{1},\ldots,q_{m})\right)  =T_{\left(
q_{1},\ldots,q_{m}\right)  }(\varphi_{z}^{-1}\circ\pi_{m})\left(
\frac{\partial}{\partial q_{a}^{k}}\right)  =(F_{k}^{a}\left(  q_{1}%
,\ldots,q_{m}\right)  ,0,\overset{a-1}{\ldots},\overset{n\text{ entries}%
}{\overbrace{(0,\overset{k-1}{\ldots},1,\ldots,0)}},\overset{m-a}{\ldots},0),
\]
where $a=1,\ldots,m$, $k=1,\ldots,n$. Observe that, by construction, the $nm$
vector fields $X_{k}^{a}$ are linearly independent and span $T_{q}%
\mathcal{G}_{z}$ at any $q\in\mathcal{G} _{z}$, since $\varphi_{z}^{-1}$ is a
diffeomorphism form $\mathbb{R}^{nm}$ to $\mathcal{G}_{z}$.

Now, we notice that for any $j=1,\ldots,l$, the vectors
\begin{equation}
\widetilde{S}_{j}\left(  X;\Gamma^{z},\Gamma_{1},\ldots,\Gamma_{m}\right)
=\left(  S_{j}\left(  X,\Gamma^{z}\right)  ,S_{j}\left(  X,\Gamma_{1}\right)
,\ldots,S_{j}\left(  X,\Gamma_{m}\right)  \right)  \label{eq 12}%
\end{equation}
are by hypothesis tangent to $\mathcal{G}_{z}$. Additionally, due to
(\ref{eq 11}) and the Stratonovich differentiation rules we can write%
\begin{equation}
\label{derivative gamma 0}\delta\Gamma_{0}=\sum_{a=1}^{m}\sum_{k=1}^{n}%
F_{k}^{a}\left(  \Gamma_{1},\ldots,\Gamma_{m}\right)  \delta\Gamma_{a}%
^{k}=\sum_{a=1}^{m}\sum_{k=1}^{n}\sum_{j=1}^{l}F_{k}^{a}\left(  \Gamma
_{1},\ldots,\Gamma_{m}\right)  S_{j}^{k}\left(  X,\Gamma_{a}\right)  \delta
X^{j}.
\end{equation}
Moreover,
\begin{equation}
\left(  \sum_{a=1}^{m}\sum_{k=1}^{n}F_{k}^{a}\left(  \Gamma_{1},\ldots
,\Gamma_{m}\right)  S_{j}^{k}\left(  X,\Gamma_{a}\right)  ,S_{j}\left(
X,\Gamma_{1}\right)  ,\ldots,S_{j}\left(  X,\Gamma_{m}\right)  \right)
\in\mathbb{R}^{n(m+1)} \label{eq 13}%
\end{equation}
belongs also to $T\mathcal{G}_{z}$ for any $j=1,\ldots,l$, since (\ref{eq 13})
can be written as a linear combination of the $nm$ linearly independent vector
fields $X_{k}^{a}$. Indeed,
\[
\left(  \sum_{a=1}^{m}\sum_{k=1}^{n}F_{k}^{a}\left(  \Gamma_{1},\ldots
,\Gamma_{m}\right)  S_{j}^{k}\left(  X,\Gamma_{a}\right)  ,S_{j}\left(
X,\Gamma_{1}\right)  ,\ldots,S_{j}\left(  X,\Gamma_{m}\right)  \right)
=\sum_{a=1}^{m}\sum_{k=1}^{n}S_{j}^{k}\left(  X,\Gamma_{a}\right)  X_{k}%
^{a}\left(  \Gamma_{1},\ldots,\Gamma_{m}\right)  .
\]
Subtracting (\ref{eq 13}) from (\ref{eq 12}), we see that for any
$j=1,\ldots,l$,
\[
W_{j}:=\left(  S_{j}\left(  X,\Gamma^{z}\right)  -\sum_{a=1}^{m}\sum_{k=1}%
^{n}F_{k}^{a}\left(  \Gamma_{1},\ldots,\Gamma_{m}\right)  S_{j}^{k}\left(
X,\Gamma_{a}\right)  ,0,\ldots,0\right)  \in T\mathcal{G}_{z}.
\]
Any of these vectors fields, if different from zero, is obviously linearly
independent from all the $X_{k}^{a}$, $a=1,\ldots,m$, $k=1,\ldots,n$. If that
is the case we could therefore conclude that $\dim(\mathcal{G}_{z})$ is
strictly bigger than $nm$, which is obviously a contradiction. Therefore,
$W_{j}=0$ necessarily, and hence
\[
S_{j}\left(  X,\Gamma^{z}\right)  =\sum_{a=1}^{m}\sum_{k=1}^{n}F_{k}%
^{a}\left(  \Gamma_{1},\ldots,\Gamma_{m}\right)  S_{j}^{k}\left(  X,\Gamma
_{a}\right)  ,
\]
which guarantees that $\Gamma_{0}$ is a solution of
(\ref{stochastic differential equation expression}) because
by~(\ref{derivative gamma 0})%
\[
\delta\Gamma_{0}=\sum_{j=1}^{l}S_{j}\left(  X,\Gamma^{z}\right)  \delta
X^{j}=\delta\Gamma^{z}.\quad\blacksquare
\]

\begin{remark}
\label{tangency remark}
\normalfont In the previous proposition we saw how the tangency of the
diagonal extensions of the vector fields that define the SDE to the
submanifolds in $\mathcal{G}$ is a sufficient condition to ensure the
existence of a superposition rule. Is it necessary? Suppose that we have a
smooth superposition rule $(\Phi,\Gamma_{1},\ldots,\Gamma_{m})$ and let $\Psi$ be
the associated map introduced in~(\ref{eq 2}). As we have that $\Psi
_{z}\left(  \Gamma^{z},\Gamma_{1},\ldots,\Gamma_{m}\right)  =0$, the
Stratonovich differentiation rules yield
\begin{equation}
0=\sum_{i=1}^{n}\sum_{a=0}^{m}\frac{\partial\Psi_{z}}{\partial q_{a}^{i}%
}\left(  \Gamma^{z},\Gamma_{1},\ldots,\Gamma_{m}\right)  \delta\Gamma_{a}%
^{i}=\sum_{j=1}^{l}\sum_{i=1}^{n}\sum_{a=0}^{m}\frac{\partial\Psi_{z}%
}{\partial q_{a}^{i}}\left(  \Gamma^{z},\Gamma_{1},\ldots,\Gamma_{m}\right)
S_{j}^{i}\left(  X,\Gamma_{a}\right)  \delta X^{j}. \label{eq 4}%
\end{equation}
A sufficient condition for this identity to hold is that, for any
$j\in\left\{  1,..,l\right\}  $,
\begin{equation}
\sum_{i=1}^{n}\sum_{a=0}^{m}\frac{\partial\Psi_{z}}{\partial q_{a}^{i}}\left(
\Gamma^{z},\Gamma_{1},\ldots,\Gamma_{m}\right)  S_{j}^{i}\left(  X,\Gamma
_{a}\right)  =0 \label{eq 5}%
\end{equation}
or, equivalently, that the diagonal extensions $\widetilde{S}_{j}\left(
X,\Gamma^{z},\Gamma_{1},\ldots,\Gamma_{m}\right)  $ are tangent to the
elements of the family of submanifolds $\mathcal{G}$ given by the zero fibers
of the maps $\Psi_{z}$. Additionally, one can find situations in
which~(\ref{eq 4}) implies~(\ref{eq 5}): for instance if $j=1$ and (like in
the case of the Brownian motion) the quadratic variation $\left[  X,X\right]
$ is a strictly increasing process, a straightforward application of the
Doob-Meyer decomposition and the It\^{o} isometry make in this
case~(\ref{eq 4}) and~(\ref{eq 5}) equivalent.
\end{remark}

\begin{remark}
\normalfont If we add to the hypotheses of Proposition~\ref{prop 2} that for any $z\in\mathbb{R}^{n}$ and for any $(p _1, \ldots ,p _m )\in \mathbb{R} ^{nm}$ there exist a submanifold $\mathcal{G}_{z}$ in $\mathcal{G}$ such that $(z,p_{1},\ldots,p_{m}%
)\in\mathcal{G}_{z}$  (for instance when $\mathcal{G}$ is a foliation of $\mathbb{R}^{n(m+1)} $ whose leaves are diffeomorphic to
$\mathbb{R}^{nm}$ via $\pi_{m}$) then the superposition function that we constructed in the proof of that result has the following extremely convenient property: the superposition function is the same for any fundamental sets of solutions $\left\{  \Gamma_{1}%
,\ldots,\Gamma_{m}\right\}  $ that we may want to choose. 
In other words, once $\Phi$ is know, we can
take $m$ arbitrary independent solutions of
(\ref{stochastic differential equation expression}) to write down any
solution. This situation  frequently occurs in mechanics; see for
instance, the study of the classical Riccati equation in (\cite{carinena-nasarre}).
\end{remark}

\section{The stochastic Lie-Scheffers Theorem}
\label{The stochastic Lie-Scheffers Theorem}

The main goal of this section is proving a theorem that characterizes the
existence of a superposition rule for a stochastic differential equation in
terms of the integrability properties of the distribution spanned by the
vector fields that define it. This can be translated into a Lie algebraic
requirement, which allows us to recover the classical Lie-Scheffers Theorem
in the stochastic context (Corollary \ref{corolario Lie clasico}).

In order to have at hand the necessary concepts to state the main theorem, we
start by briefly recalling some standard results on generalized distributions
due to Stefan~\cite{Stefan1974a, Stefan1974b} and Sussman~\cite{Sussman1973}.
Let $M$ be a smooth manifold, $\mathcal{D}\subset\mathfrak{X}(M)$ be a family
of smooth vector fields, and $D$ the smooth generalized distribution spanned
by $\mathcal{D}$. Let $G_{\mathcal{D}}$ be the pseudogroup of transformations
generated by the flows of the vector fields in $\mathcal{D}$ and constructed
as follows: let $k\in\mathbb{N}^{\ast}$ be a positive natural number,
$\mathcal{X}$ an ordered family $\mathcal{X}=(X_{1},\ldots,\,X_{k})$ of $k$
elements of $\mathcal{D}$, and $T$ a $k$--tuple $T=(t_{1},\ldots,\,t_{k}%
)\in\mathbb{R}^{k}$ such that $F_{t}^{i}$ denotes the (locally defined) flow
of $X_{i}$, $i\in\{1,\ldots,\,k\}$, $t_{i}$; the elements $\mathcal{F}_{T}$ of
$G_{\mathcal{D}}$ are the locally defined diffeomorphisms of the form
$\mathcal{F}_{T}=F_{t_{1}}^{1}\circ F_{t_{2}}^{2}\circ\cdots\circ F_{t_{k}%
}^{k}$. Two points $x$ and $y$ in $M$ are said to be $G_{\mathcal{D}}%
$-equivalent, if there exists a diffeomorphism $\mathcal{F}_{T}\in
G_{\mathcal{D}}$ such that $\mathcal{F}_{T}(x)=y$. The relation
$G_{\mathcal{D}}$--equivalent is an equivalence relation whose equivalence
classes are called the $G_{\mathcal{D}}$-orbits, that are sometimes referred to as the {\bfi  accessible sets} associated to the family $\mathcal{D}$.

Given the family $\mathcal{D}$ and the associated pseudogroup $G_{\mathcal{D}%
}$ we can define another family $\mathcal{D}^{\prime}$ of vector fields as
\[
\mathcal{D}^{\prime}:=\{T\mathcal{F}_{T}\cdot X\mid X\in\mathcal{D}%
,\mathcal{F}_{T}\in G_{\mathcal{D}}\},
\]
that clearly extends $\mathcal{D}$, that is, $\mathcal{D}\subset
\mathcal{D}^{\prime}$. The distribution $D^{\prime}$ spanned by the elements
of $\mathcal{D}^{\prime}$ is by construction $G_{\mathcal{D}}$-invariant. That
is, for each $\mathcal{F}_{T}\in G_{\mathcal{D}}$ and for each $z\in M$ in the
domain of $\mathcal{F}_{T}$,%
\begin{equation}
T_{z}\mathcal{F}_{T}(D^{\prime}(z))=D^{\prime}(\mathcal{F}_{T}(z)).
\label{eq 19}%
\end{equation}
Moreover, since $(\mathcal{D}^{\prime})^{\prime}=\mathcal{D}^{\prime}$ by
construction, the Stefan-Sussmann Theorem guarantees that it is {\bfi
completely integrable} in the sense that for every point $z\in M$, there
exists an integral manifold of $D^{\prime}$ everywhere of maximal dimension
which contains $z$. The maximal integral manifolds of a completely integrable
generalized distribution on $M$ form a {\bfi generalized foliation} of $M$
(see for instance~\cite{dazord 1985}). A leaf of a generalized foliation is
{\bfi regular} if it has a neighborhood where the singular foliation induces a
regular foliation by restriction. A point is regular if it belongs to a
regular leaf. Regular points are open and dense in $M$ (\cite[Th\'eor\`eme
2.2]{dazord 1985}). We will refer to $D^{\prime}$ (respectively $\mathcal{D}%
^{\prime}$) as the {\bfi Stefan-Sussmann extension} of $D$ (respectively
$\mathcal{D}$). The Stefan-Sussmann's Theorem also establishes an equivalence
between the $G_{\mathcal{D}}$-invariance of $D$ ($D^{\prime}=D$) and its
complete integrability; additionally, if $D$ is a completely integrable
distribution, then its integral manifolds are the $G_{D}$-orbits. When the distribution $D$ has constant dimension, the Stefan-Sussmann Theorem reduces to the celebrated and especially convenient Frobenius Theorem which states the $D$ is integrable if and only if $D$ is
involutive. Recall that $D$ is {\bfi involutive} if
$[X,\,Y]$ takes values in $D$ whenever $X$ and $Y$ are vector fields with
values in $D$.

\medskip

In the sequel, we will use the following notation in order to be able to handle
diagonal extensions of different dimensions. Given $l\in\mathbb{N}$ and
$X\in\mathfrak{X}(\mathbb{R}^{n})$, we will denote by $\widetilde{X}^{l}%
\in\mathfrak{X}(\mathbb{R}^{ln})$ the diagonal extension of $X$ to
$\mathbb{R}^{ln}$. For the sake of consistency with the previous section
$\widetilde{X}$ means $\widetilde{X}^{m+1}$.

\begin{theorem}
[Lie-Scheffers' Theorem for SDE]
\label{teorema Lie}
Let
\begin{equation}
\delta\Gamma=S\left(  X,\Gamma\right)  \delta X \label{eq 1bis}%
\end{equation}
be a stochastic differential equation on $\mathbb{R}^{n}$, where
$X:\mathbb{R}_{+}\times\Omega\rightarrow\mathbb{R}^{l}$ is a given
$\mathbb{R}^{l}$-valued semimartingale and $S\left(  x,z\right)
:T_{x}\mathbb{R}^{l}\longrightarrow T_{p}\mathbb{R}^{n}$ is a Stratonovich
operator from $\mathbb{R}^{l}$ to $\mathbb{R}^{n}$. Let $V$ be an arbitrary open
neighborhood of $\mathbb{R}^{n}$. Then,

\begin{enumerate}
\item[{\bf (i)}] If the $X$-dependent vector fields $\left\{  S_{1}\left(
X,\cdot\right)  ,\ldots,S_{l}\left(  X,\cdot\right)  \right\}  $ can be
expressed on $V$ as
\begin{equation}
S_{j}\left(  X,z\right)  =\sum_{i=1}^{r}b_{j}^{i}\left(  X\right)
Y_{i}\left(  z\right)  \in T_{z}\mathbb{R}^{n},\quad b_{j}^{i}\in C^{\infty
}(\mathbb{R}^{l}),\quad z\in V, \label{eq 7}%
\end{equation}
and the distribution $D$ spanned by the vector fields $\mathcal{D}=\{Y_{1},\ldots,Y_{r}\}\subset
\mathfrak{X}\left(  V\right)  $ is involutive, then (\ref{eq 1bis}) admits a local
superposition rule.

\item[{\bf (ii)}] Conversely, suppose that (\ref{eq 1bis}) admits a superposition
rule $(\Phi,\{\Gamma_{1},\ldots,\Gamma_{m}\})$ and that the diagonal
extensions $\{\widetilde{S}_{1}\left(  X,\cdot\right)  ,\allowbreak
\ldots,\widetilde{S}_{l}\left(  X,\cdot\right)  \}$ to $\mathbb{R}^{n(m+1)}$
are tangent to the family $\mathcal{G}$ of $nm$-dimensional submanifolds of
$\mathbb{R}^{n(m+1)}$ associated to this superposition rule (see
Proposition~\ref{prop 2}). Let $\widetilde{D}(q):=\mathrm{span}\{\widetilde
{S}_{j}(X_{t},q)\mid j\in\{1,\ldots,l\},\,t\in\mathbb{R}_{+}\}$,
$q\in\mathbb{R}^{n(m+1)}$, $\widetilde{D}^{\prime}$ the Stefan-Sussmann
extension of $\widetilde{D}$, and $\mathcal{G}_{0}$ its associated generalized
foliation. Let $z\in\mathbb{R}^{n}$, $p_{i}=(\Gamma_{i})_{t=0}$, and suppose
that $p=(z,p_{1},\ldots,p_{m})\in\mathbb{R}^{n(m+1)}$ belongs to a regular
leaf $\left(  \mathcal{G}_{0}\right)  _{z}$ of $\mathcal{G}_{0}$. Then, there
exists an open neighborhood $V$ of $z$, a family of vector fields $\left\{
Y_{1},\ldots,Y_{r}\right\}  \subset\mathfrak{X}\left(  V\right)  $, and a
family of functions $\{b_{j}^{i}\}_{j=1,..,l}^{i=1,..,r}\subset C^{\infty
}\left(  \mathbb{R}^{l}\right)  $ such that%
\begin{equation}
\label{decomposition with the Ys}
S_{j}\left(  X,v\right)  =\sum_{i=1}^{r}b_{j}^{i}\left(  X\right)
Y_{i}\left(  v\right),  
\end{equation}
for any $v\in V$. Moreover, the vector fields $\left\{  Y_{1},\ldots
,Y_{r}\right\}  $ form a real Lie algebra.
\end{enumerate}
\end{theorem}

\noindent\textbf{Proof.\ \ (i)} Given $l\in\mathbb{N}$, we define
$V^{l}:=V\times\overset{l)}{\ldots}\times V$ and $d_{l}:=\max_{q\in V^{l}%
}\left\{  \dim\left(  \mathrm{span}\{\widetilde{Y}_{1}^{l}(q),\ldots
,\widetilde{Y}_{r}^{l}(q)\}\right)  \right\}  $. Notice that for any
$l\in\mathbb{N}$ one has $d_{l}\leq d_{l+1}$ and $d_{l}\leq r$. Let
$m\in\mathbb{N}$ be the smallest number for which $d_{m}=d_{m+1}$ and let
$q_{0}\in V^{m+1}$ be such that
\begin{equation}
\dim\left(  \mathrm{span}\{\widetilde{Y}_{1}^{m+1}(q_{0}),\ldots,\widetilde
{Y}_{r}^{m+1}(q_{0})\}\right)  =d_{m+1}. \label{definition z 0}%
\end{equation}
The maximality of the dimension of $\mathrm{span}\{\widetilde{Y}_{1}%
^{m+1},\ldots,\widetilde{Y}_{r}^{m+1}\}$ at $q_{0}$ implies that there exists
a neighborhood $U$ of $q_{0}$ in $V^{m+1}$ for which $\dim(\mathrm{span}%
\{\widetilde{Y}_{1}^{m+1}(q),\ldots,\widetilde{Y}_{r}^{m+1}(q)\})=d_{m+1}$,
for all $q\in U$. Indeed, the expression~(\ref{definition z 0}) is equivalent
to saying that the $r\times n(m+1)$ matrix $M(q)$ with entries $M_{ij}%
(q):=(\widetilde{Y}_{i}^{m+1}(q))^{j}$ has rank $d_{m}$ when evaluated at
$q_{0}$ which, in turn, amounts to the existence of a non-vanishing minor
$M_{d_{m+1}}(q_{0})$ of $M(q_{0})$ of order $d_{m+1}$. Since the minor
$M_{d_{m+1}}(q)$ depends smoothly on $q$ and $M_{d_{m+1}}(q_{0})\neq0$, there
exists an open neighborhood $U$ of $q_{0}$ in $V^{m+1}$ for which $M_{d_{m+1}%
}(q)\neq0$, for any $q\in U$. This implies that $\dim(\mathrm{span}%
\{\widetilde{Y}_{1}^{m+1}(q),\ldots,\widetilde{Y}_{r}^{m+1}(q)\})\geq d_{m+1}%
$, for all $q\in U$. However, the maximality used in the definition of
$d_{l+1}$ implies that the previous inequality is necessarily an equality.

Consequently, we have found an open set $U\subset V^{m+1}$ in which the
distribution $D$ spanned by the family $\{\widetilde{Y}_{1}^{m+1}%
,\ldots,\widetilde{Y}_{r}^{m+1}\}$ has constant rank. Moreover, (\ref{eq 14})
and the hypothesis on $\{Y_{1},\ldots,Y_{r}\}$ being in involution imply by
the classical Frobenius Theorem that $D$ is integrable. Let $\mathcal{G}_{0}$
be the family of maximal integrable leaves of $D$ that form a foliation of
$U^{m+1}$. Now, shrinking $U$ if necessary and using foliation coordinates for
$\mathcal{G}_{0}$, we extend the distribution $D$ to another integrable
distribution $\overline{D}\supset D$ of rank $nm$ whose integrable leaves
$\mathcal{G}$ contain those of $\mathcal{G}_{0}$, and for which the
restrictions of $\pi_{m}:\mathbb{R}^{n(m+1)}\rightarrow\mathbb{R}^{mn}$ to the
leaves in $\mathcal{G}$ are diffeomorphisms onto their images.

Let now $\{p_{1},\ldots,p_{m}\}$ be a set of $m$ distinct points in $V$ such
that $(p_{1},\ldots,p_{m})\in\pi_{m}(U)$ and $\left\{  \Gamma_{1}%
,\ldots,\Gamma_{m}\right\}  $ the solutions of of (\ref{eq 1bis}) such that
$(\Gamma_{1},\ldots,\Gamma_{m})_{t=0}=(p_{1},\ldots,p_{m})$ a.s.. Let
$\Gamma:=(\Gamma_{1},\ldots,\Gamma_{m})$ and $\tau$ the stopping time defined
as $\tau:=\inf\{t>0\mid\Gamma_{t}\neq\pi_{m}(U)\}$. Since the vector fields
\[
\widetilde{S}_{j}^{m+1}\left(  X,\Gamma\right)  =\sum_{i=1}^{r}b_{j}%
^{i}\left(  X\right)  \widetilde{Y}_{i}^{m+1}\left(  \Gamma\right)
\]
are tangent to the integral leaves of $\mathcal{G}_{0}$ and hence to those of
$\mathcal{G}$, at least up to time $\tau$, Proposition~\ref{prop 2} guarantees
the existence of a local superposition rule.

\medskip

\noindent\textbf{(ii)} We start the proof by providing a lemma that will be
needed in our argument.

\begin{lemma}
\label{lemma 1} Let $\{Y_{1},\ldots,Y_{r}\}\subset\mathfrak{X}(\mathbb{R}%
^{n})$ with $r\leq mn$ and let $\{\widetilde{Y}_{1},\ldots,\widetilde{Y}%
_{r}\}$ be the corresponding diagonal extensions to $\mathbb{R}^{n(m+1)}$.
Suppose that $\{T_{q}\pi_{m}(\widetilde{Y}_{1}(q)),\ldots,T_{q}\pi
_{m}(\widetilde{Y}_{r}(q))\}$ are linearly independent for any $q$ in a
neighborhood $U\subseteq\mathbb{R}^{n(m+1)}$. If the sum $\sum_{i=1}^{r}%
b^{i}\widetilde{Y}_{i}$ with $b^{i}\in C^{\infty}\left(  U\right)  $,
$i=1,\ldots,r$, is again a diagonal extension then the functions $b^{i}$ are
necessarily the pull-back by $\pi_{m}$ of a family functions in $C^{\infty
}(\pi_{m}\left(  U\right)  )$. More specifically, if $(q_{a}^{j}%
;j=1,\ldots,n;a=0,\ldots,m)$ are coordinates for $\mathbb{R}^{n(m+1)}$, then
the functions $\{b^{i}\}_{i=1,..,r}$ do not depend on $(q_{0}^{j}%
;j=1,\ldots,n)$.
\end{lemma}

\noindent\textbf{Proof.\ \ } Using the coordinates $\left(  q^{j}%
;j=1,\ldots,n\right)  $ for $\mathbb{R}^{n}$, there exists a family of
functions $A_{i}^{j}\in C^{\infty}(\mathbb{R}^{n})$, $i\in\{1,\ldots,r\}$ ,
$j\in\{1,\ldots,n\}$, such that the vector fields $\{Y_{1},\ldots
,Y_{r}\}\subset\mathfrak{X}(\mathbb{R}^{n})$ can be written as%
\[
Y_{i}(q)=\sum_{j=1}^{n}A_{i}^{j}\left(  q\right)  \frac{\partial}{\partial
q^{j}}%
\]
which implies that the diagonal extensions have the expression
\[
\widetilde{Y}_{i}(q_{0},\ldots,q_{m})=\sum_{a=0}^{m}\sum_{j=1}^{n}A_{i}%
^{j}\left(  q_{a}\right)  \frac{\partial}{\partial q_{a}^{j}}.
\]
Then, if we assume that
\[
\sum_{i=1}^{r}b^{i}\left(  q_{0},\ldots,q_{m}\right)  \widetilde{Y}_{i}%
(q_{0},\ldots,q_{m})=\sum_{i=1}^{r}\sum_{a=0}^{m}\sum_{j=1}^{n}b^{i}\left(
q_{0},\ldots,q_{m}\right)  A_{i}^{j}\left(  q_{a}\right)  \frac{\partial
}{\partial q_{a}^{j}}%
\]
is a diagonal extension on $U$, then there exist some functions $\{B^{i}%
\}_{i=1,\ldots,r}\subset C^{\infty}\left(  \mathbb{R}^{n}\right)  $ such that
\[
\sum_{i=1}^{r}\left.  b^{i}\left(  q_{0},\ldots,q_{m}\right)  A_{i}^{j}\left(
q_{a}\right)  \right\vert _{U}=\left.  B^{j}(q_{a})\right\vert _{U},\text{
\ \ }a=0,\ldots,m,~~j=1,\ldots,n.
\]
That is, the $r$ functions $b^{i}\left(  q_{0},\ldots,q_{m}\right)  $ solve
the following subsystem of linear equations
\begin{equation}
\left(
\begin{array}
[c]{c}%
\mathcal{A}(q_{0})\\
\mathcal{A}(q_{1})\\
\vdots\\
\mathcal{A}(q_{m})
\end{array}
\right)  \left(
\begin{array}
[c]{c}%
b^{1}(q_{0},\ldots,q_{m})\\
\vdots\\
b^{r}(q_{0},\ldots,q_{m})
\end{array}
\right)  =\left(
\begin{array}
[c]{c}%
\mathcal{B}(q_{0})\\
\mathcal{B}(q_{1})\\
\vdots\\
\mathcal{B}(q_{m})
\end{array}
\right)  \label{system linear 1}%
\end{equation}
where $\mathcal{A}$ and $\mathcal{B}$ are the $n(m+1)\times r$ and
$n(m+1)\times1$ matrices, respectively, defined as $\mathcal{A}(q_{a}%
)_{ij}=A_{j}^{i}(q_{a})$ and $\mathcal{B}(q_{a})_{i}=B^{i}(q_{a})$,
$a=0,\ldots,m$. Now, the hypothesis on the linear independence of $\{T\pi
_{m}(\widetilde{Y}_{1}),\ldots,T\pi_{m}(\widetilde{Y}_{r})\}$ implies that the
rank of the matrix $\left(  \mathcal{A}(q_{1}),\ldots,\mathcal{A}%
(q_{m})\right)  $ is $r\leq nm$ and hence~(\ref{system linear 1}) has a unique
solution which coincides with the unique solution of the system
\begin{equation}
\left(
\begin{array}
[c]{c}%
\mathcal{A}(q_{1})\\
\vdots\\
\mathcal{A}(q_{m})
\end{array}
\right)  \left(
\begin{array}
[c]{c}%
b^{1}(q_{0},\ldots,q_{m})\\
\vdots\\
b^{r}(q_{0},\ldots,q_{m})
\end{array}
\right)  =\left(
\begin{array}
[c]{c}%
\mathcal{B}(q_{1})\\
\vdots\\
\mathcal{B}(q_{m})
\end{array}
\right)  . \label{system linear 2}%
\end{equation}
Since there is no dependence on the coordinates $q_{0}$ in the augmented
matrix associated to the system~(\ref{system linear 2}), its solution
$(b^{1},\ldots,b^{r})$ does not therefore depend on $q_{0}$, as required.
\quad$\blacktriangledown$

\medskip

Suppose now that the stochastic differential equation~(\ref{eq 1bis}) admits a
superposition rule and that we are in the hypotheses of the theorem. We start
by emphasizing that since the vector fields $\{\widetilde{S}_{1}\left(
X,\cdot\right)  ,\allowbreak\ldots,\widetilde{S}_{l}\left(  X,\cdot\right)
\}$ are, by hypothesis, tangent to the elements of the family $\mathcal{G}$
then their flows leave invariant those submanifolds and hence, the
Stefan-Sussmann extension $\widetilde{D}^{\prime}$ of $\widetilde{D}$ is also
tangent to the elements of $\mathcal{G}$. This argument guarantees that, given
the regular leaf $(\mathcal{G}_{0})_{z}$ of $\mathcal{G}_{0}$, then there
exists an element $\mathcal{G}_{z}$ in $\mathcal{G}$ that contains it.

Now since $p=(z,p_{1},\ldots,p_{m})\in\mathbb{R}^{n(m+1)}$ belongs to a
regular leaf $\left(  \mathcal{G}_{0}\right)  _{z}$ of $\mathcal{G}_{0}$, then
there is an open neighborhood $U$ of $p$ where we can choose (taking regular
foliation coordinates) a family of linearly independent vector fields
$\{\widetilde{Y}_{1},\ldots,\widetilde{Y}_{r}\}\subset\mathfrak{X}%
(\mathbb{R}^{n(m+1)})$ that span the tangent spaces to the leaves of
$\mathcal{G}_{0}\cap U$. The vector fields $\{\widetilde{Y}_{1},\ldots
,\widetilde{Y}_{r}\}$ can be chosen as the diagonal extensions of $r$ vector
fields $\{Y_{1},\ldots,Y_{r}\}\subset\mathfrak{X}(\mathbb{R}^{n})$, since the
Stefan-Sussmann extension $\widetilde{D}^{\prime}=\mathrm{span}\{T\widetilde
{\mathcal{F}}_{T}\cdot\widetilde{S}_{i}\left(  X,\cdot\right)  \mid
i\in\{1,\ldots,l\},\widetilde{\mathcal{F}}_{T}\in G_{\mathcal{D}}\}$ of
$\widetilde{D}$ is made of diagonal extensions. Indeed, in order to see that
$\widetilde{D}^{\prime}$ is spanned by diagonal extensions, it suffices to
notice that the flow $\widetilde{F}_{t}$ of the diagonal extension
$\widetilde{Y}\in\mathfrak{X}(\mathbb{R}^{n(m+1)})$ of a vector field
$Y\in\mathfrak{X}(\mathbb{R}^{n})$ is $\widetilde{F}_{t}(q_{0},\ldots
,q_{m})=(F_{t}(q_{0}),\ldots,F_{t}(q_{m}))$, with $F_{t}$ the flow of $Y$;
hence%
\begin{align*}
T_{q}\widetilde{F}_{t}(\widetilde{Y}(q))  &  =\left(T_{q_{0}}F_{t}\times\ldots\times
T_{q_{m}}F_{t}\right)(\widetilde{Y}\left(  q\right)  )\\
&  =(T_{q_{0}}F_{t}(Y(q_{0})),\ldots,T_{q_{m}}F_{t}(Y\left(  q_{m}\right)
)=\widetilde{\left(  TF_{t}(Y)\right)  }(q)
\end{align*}
is again a diagonal extension. Given that by (\ref{eq 14}) diagonal extensions
form an algebra, the statement follows.

Moreover, since the distribution $\widetilde{D}^{\prime}|_{U}$ is regular and
integrable then it is necessarily integrable in the sense of Frobenius, that
is, there exist functions $\{c_{ij}^{k}\}_{i,j,k=1,..,r}\subset C^{\infty
}(\mathbb{R}^{n(m+1)})$ such that
\begin{equation}
\left[  \widetilde{Y}_{j},\widetilde{Y}_{i}\right]  =\sum_{k=1}^{r}c_{ji}%
^{k}\widetilde{Y}_{k}. \label{eq 20}%
\end{equation}
Now, as $[\widetilde{Y}_{j},\widetilde{Y}_{i}]=\widetilde{[Y_{j},Y_{i}]}$, we
conclude that $\sum_{k=1}^{r}c_{ji}^{k}\widetilde{Y}_{k}$ is a diagonal
extension. Also, as the projection $\pi_{m}$ is a local diffeomorphism when
restricted to $U\cap\mathcal{G}_{z}$, the family of vectors $\{T\pi
_{m}(\widetilde{Y}_{1}),\ldots,T\pi_{m}(\widetilde{Y}_{r})\}$ is necessarily
linearly independent. In these circumstances Lemma \ref{lemma 1} implies that
the coefficients $\{c_{ij}^{k}\}_{i,j,k=1,..,r}$ do not depend on the first
$n$ coordinates $q_{0}^{j}$, $j=1,\ldots,n$. We now apply $\pi_{\mathbb{R}%
^{n}}^{0}$ (see (\ref{eq 15})) on both sides of~(\ref{eq 20}) and we obtain%
\begin{equation}
\left[  Y_{j},Y_{i}\right]  (v)=\sum_{k=1}^{r}c_{ji}^{k}(q_{1},\ldots
,q_{m})Y_{k}(v) \label{eq 20 bis}%
\end{equation}
where $v\in V:=\pi_{\mathbb{R}^{n}}^{0}(U)$ and $\left(  q_{1},\ldots
,q_{m}\right)  \in\mathbb{R}^{nm}$ is any arbitrary point such that $\left(
v,q_{1},\ldots,q_{m}\right)  \in U$. Since the left hand side of
(\ref{eq 20 bis}) does not depend on $\left(  q_{1},\ldots,q_{m}\right)  $ then the dependence of the coefficients $c_{ji}^{k}(q_{1},\ldots,q_{m})$ on those coordinates is necessarily trivial
which allows us to
conclude that $\{Y_{1},\ldots,Y_{r}\}$ close a Lie algebra.

Finally, since the vector fields $\widetilde{S}_{j}\left(  X,\cdot\right)  $
are tangent to $\mathcal{G}_{0}$, $j=1,\ldots,l$, then there is a family of
$X$-dependent functions $b_{j}^{i}\left(  X,\cdot\right)  \in C^{\infty
}\left(  U\right)  $ such that%
\[
\widetilde{S}_{j}\left(  X,q\right)  =\sum_{i=1}^{r}b_{j}^{i}\left(
X,q\right)  \widetilde{Y}_{i}\left(  q\right)  ,
\]
for any $q\in U$. As $\widetilde{S}_{j}\left(  X,\cdot\right)  $ is also a
diagonal extension, we can use again Lemma~\ref{lemma 1} in order to prove
that the functions $\{b_{j}^{i}\}_{j=1,..,l}^{i=1,..,r}$ do not depend on
$q_{0}$. Consequently,
\begin{equation}
\widetilde{S}_{j}\left(  X,q\right)  =\sum_{i=1}^{r}b_{j}^{i}\left(
X,(q_{1},\ldots,q_{m})\right)  \widetilde{Y}_{i}\left(  p\right)  .
\label{eq 16}%
\end{equation}
As we did in the previous paragraph, we  apply $\pi_{\mathbb{R}^{n}}^{0}$ on both sides of~(\ref{eq 16}) 
\[
S_{j}\left(  X,v\right)  =\sum_{i=1}^{r}b_{j}^{i}\left(  X,(q_{1},\ldots
,q_{m})\right)  Y_{i}\left(  v\right)  ,
\]
for any $v\in V$. Again, we realize that since the left hand side of this equation is independent of
$(q_{1},\ldots,q_{m})$, the dependence of the functions $b_{j}^{i}$ on the coordinates $\left(  q_{1},\ldots,q_{m}\right)  $ is necessarily trivial, which yields expression~(\ref{decomposition with the Ys}). \quad$\blacksquare$

\bigskip

\begin{remark}
\normalfont Theorem \ref{teorema Lie} is a generalization for stochastic
differential equations of the classical Lie-Scheffers  Theorem stated for
time-dependent ordinary differential equations. That theorem claims that a
differential equation $\dot{y}=Y\left(  t,y\right)  $
on $\mathbb{R}^{n}$ given by a time-dependent vector field $Y\left(
t,\cdot\right)  \in\mathfrak{X}\left(  \mathbb{R}^{n}\right)  $,
$t\in\mathbb{R}$, admits a superposition rule if and only if $Y$ can be
locally written in the form $Y\left(  t,y\right)  =\sum_{i=1}^{r}f^{i}%
(t)Y_{i}(y)$, where $\{f^{i}\}_{i=1,\ldots,r}\subset C^{\infty}(\mathbb{R})$
and $\{Y_{1},\ldots,Y_{r}\}\subset\mathfrak{X}\left(  \mathbb{R}^{n}\right)  $
form a (real) Lie subalgebra of $(\mathfrak{X}(M),[\cdot,\cdot])$ (see
\cite{Cari07} and \cite{libro marmo}). In relation to the traditional
presentation of the Lie-Scheffers Theorem, our Theorem \ref{teorema Lie}:

\begin{enumerate}
\item[{\bf (i)}] weakens the hypotheses under which we can guarantee the existence
of superposition rules. 
The involutivity of the vector fields $\{Y_{1}%
,\ldots,Y_{r}\}$ is, in general, less restrictive than requiring that they form a 
Lie algebra over the reals. We know \emph{a posteriori} by the second part of Theorem \ref{teorema Lie}
that, around regular points, if there exists a superpositon rule, the 
components $\{S_{1},\ldots,S_{l}\}$ of the Stratonovich operator can also be expressed in terms of a family of vector fields that close  a Lie algebra.

\item[{\bf (ii)}] carefully spells out the regularity conditions under which we have a converse; those conditions are only vaguely evoked in the already cited deterministic papers.

\item [{\bf (iii)}] It is worth noticing that, apart from the two points that we just explained,~Theorem \ref{teorema Lie} cannot be seen as a mere transcription of the
deterministic Lie-Scheffers Theorem into the context of Stratonovich stochastic integration by using the so called  Malliavin's Transfer
Principle~\cite{malliavin transfer} due to the purely stochastic conditions that appear in the statement of the theorem. Those additional requirements have to do with the tangency of the diagonal extensions of the components of the Stratonovich operator to the family of submanifolds associated to the superposition rule (see also Remark~\ref{tangency remark}).
\end{enumerate}
\end{remark}

In the next corollary, we show for the sake of completeness how the classical statement of the
Lie-Scheffers Theorem (generalized to SDEs) can be easily obtained out of Theorem~\ref{teorema Lie}.

\begin{corollary}
\label{corolario Lie clasico}
Using the notation in Theorem
\ref{teorema Lie}, suppose that the $X$-dependent family of vector fields $\left\{
S_{1}\left(  X,\cdot\right)  ,\ldots,S_{l}\left(  X,\cdot\right)  \right\}  $
that define the stochastic differential equation
(\ref{stochastic differential equation expression}) can be expressed as
\[
S_{j}\left(  X,z\right)  =\sum_{i=1}^{r}b_{j}^{i}\left(  X\right)
Y_{i}\left(  z\right)  \in T_{z}\mathbb{R}^{n},\quad b_{j}^{i}\in C^{\infty
}(\mathbb{R}^{l}),\quad z\in \mathbb{R}^n.
\]
Let $\mathrm{Lie}\{Y_{1},\ldots,Y_{r}\}$ be the real Lie subalgebra of
$(\mathfrak{X}(\mathbb{R}^{n}),[\cdot,\cdot])$ generated by the family
$\{Y_{1},\ldots,Y_{r}\}\subset\mathfrak{X}(\mathbb{R} ^n)$. If $\mathrm{Lie}
\{Y_{1},\ldots,Y_{r}\}$ is finite dimensional then
(\ref{stochastic differential equation expression}) has a superposition rule.
\end{corollary}

\noindent\textbf{Proof.\ \ } Let $D$ and $D_{2}$ be the generalized distributions
associated to the families of vector fields $\mathcal{D}=\{Y_{1},\ldots,Y_{r}\}$ and $\mathcal{D}%
_{2}=\mathrm{Lie}\{Y_{1},\ldots,Y_{r}\}$, respectively. Observe that if
$D(z)\varsubsetneq D_{2}(z)$, $z\in \mathbb{R}^n$, then since $\mathrm{Lie}%
\{Y_{1},\ldots,Y_{r}\}$ is finite dimensional, we can always complete the family
$\{Y_{1},\ldots,Y_{r}\}$ with a finite number of vectors $\{Z_{1},\ldots
,Z_{s}\}\subset\mathcal{D}$ such that $D(z)=D_{2}(z)$. We then write the $X$-dependent vector
fields $\{S_{1}\left(  X,\cdot\right)  ,\ldots,S_{l}\left(  X,\cdot\right)
\}$ as%
\[
S_{j}\left(  X,z\right)  =\sum_{i=1}^{r}b_{j}^{i}\left(  X\right)
Y_{i}\left(  z\right)  +\sum_{k=1}^{s}a_{j}^{k}(X)Z_{k}(z),\quad z\in \mathbb{R}^n,
\]
with $a_{j}^{k}=0$ for any $j=1,\ldots,l$ and any $k=1,\ldots,s$.
Therefore, we may simply suppose that $D(z)=\operatorname{span}\{\mathrm{Lie}
\{Y_{1},\ldots,Y_{r}\}(z)\}$, $z\in \mathbb{R}^n$ and since $D_{2}$ is trivially involutive,
the corollary follows from Theorem \ref{teorema Lie} (i). \quad$\blacksquare$

\section{Lie-Scheffers systems and stochastic differential equations on Lie groups and homogeneous spaces}

\label{Lie-Scheffers systems on Lie groups and homogeneous spaces}

The Lie-Scheffers systems that are defined by a set of vector fields that generate a finite dimensional Lie algebra, that is, those that satisfy the hypothesis of
Corollary \ref{corolario Lie clasico} or of Theorem \ref{teorema Ben Arous} can be reformulated in the language of group actions. More specifically, as we see in the next proposition, such systems come down locally to studying the solutions of an equivalent Lie-Scheffers system on a Lie group.

\begin{proposition}
\label{locally on lie groups}
Consider a stochastic differential equation that satisfies the hypotheses of Corollary \ref{corolario Lie clasico}. Let $z\in M$ be a point such that there exists a neighborhood $V$ of $z$ in which the dimension of $ {\rm Lie} \left\{Y _1, \ldots, Y _r \right\}$ is constant. Then, shrinking $V$ if necessary, there exists a Lie group $G$ such that  $\dim\left(
G\right)  =\dim\left(  \mathrm{Lie}\{Y_{1},\ldots,Y_{r}\}|_V\right)  $, a group action $ \Xi:G \times  V \rightarrow V $, and Lie algebra elements $\{\xi_{1},\ldots,\xi_{r}\}\subset\mathfrak{g}$ such that
\begin{equation}
\label{action and infinitesimal generators}
Y_{i}(z)=\xi_{i}^{M}(z):=\left.  \frac{d}{dt}\right\vert _{t=0}\Xi\left(
\exp\left(  t\xi_{i}\right)  ,z\right)  ,~~z\in V.
\end{equation}
Moreover, the solution  starting at $z\in M$ of the restriction to $V$ of the stochastic differential equation may be expressed as
\begin{equation}
\Gamma_{t}^{z}=\Xi\left(  g_{t},z\right)  , \label{eq 38}%
\end{equation}
where $g_{t}:\mathbb{R}_{+}\times\Omega\rightarrow G$ is the semimartingale
solution of the stochastic differential equation on $G$%
\begin{equation}
\delta g_{t}=\sum_{i=1}^{r}\xi_{i}^{G}\left(  g_t\right)  \delta X_{t}^{i}
\label{eq 40}%
\end{equation}
with initial condition $g_{t=0}=e$ a.s.
\end{proposition}

\noindent\textbf{Proof.\ \ }
Since the statement of the proposition is local 
we can always assume that the vector fields $\{Y_{1}
,\ldots,Y_{r}\}$ are complete by multiplying them by a compactly supported bump function and by restricting ourselves to an open neighborhood $V$ consistent with that construction. In that situation and if $\dim\left(  \mathrm{Lie}\{Y_{1},\ldots,Y_{r}\}|_V\right)  < \infty$, Palais showed in
\cite{Palais} (see Corollary in page 97 and Theorem III in page 95) that there exists a unique connected Lie group $G$ contained in
the group of diffeomorphisms of $M$ and a left action $\Xi:G\times
M\rightarrow M$ such that~(\ref{action and infinitesimal generators}) holds and $T_{e}\Xi_{z}:\mathfrak{g}\rightarrow\mathrm{Lie}%
\{Y_{1},\ldots,Y_{r}\}(z)$ is an isomorphism, for any $z\in V$. 

Let now $g_{t}:\mathbb{R}_{+}\times\Omega\rightarrow G$ be the solution semimartingale
of the stochastic differential equation on $G$
\begin{equation}
\delta g_{t}=\sum_{i=1}^{r}\xi_{i}^{G}\left(  g_t\right)  \delta X_{t}^{i},
\label{eq 40 star}%
\end{equation}
where $\xi_{i}^{G}\in\mathfrak{X}
\left(  G\right)  $ denotes the right invariant infinitesimal generator associated to
$\xi_{i}\in\mathfrak{g}$ via the left translations of $G$ on $G$. Given that any two infinitesimal generators $\xi^{G}$ and $\xi^{M}$, $\xi\in \mathfrak{g}$, are related by the formula $T_{g}\Xi_{z}(\xi^{G}
)=\xi^{M}(\Xi\left(  g,z\right)  )$, $g \in  G $, $z \in V $, it is straightforward to verify that if  $g_{t}$ is a solution of (\ref{eq 40}) with initial condition $g_{t=0}=e$ a.s., then 
\begin{equation*}
\Gamma_{t}^{z}=\Xi\left(  g_{t},z\right)  ,
\end{equation*}
is the solution of $
\delta\Gamma_{t}=\sum_{i=1}^{r}Y_{i}\left(  \Gamma_{t}\right)  \delta
X_{t}^{i}$ such that  $ \Gamma_0=z $, a.s.
\quad $\blacksquare$

\begin{remark}
\normalfont
Observe that (\ref{eq 38}) may be
understood as a general reformulation of (\ref{eq 43}) (see also
\cite[Th\'{e}or\`{e}me 19]{Ben Arous}). Processes of the type
$\Gamma_{t}^{z}=\Xi\left(  g_{t},z\right)  $ defined using a group action
are sometimes called {\bfi one point motions} (\cite{Liao book}). 
\end{remark}

The proposition that we just proved shows that for Lie-Scheffers systems defined by vector fields that generate a finite dimensional Lie algebra $\mathfrak{g}$, it is the 
associated Lie-Scheffers system  on the Lie group $G$ (\ref{eq 40}) that really
matters. This is the subject of the rest of this section.

\medskip

\noindent {\bf Stochastic differential equations on Lie groups.} Let now $G$ be an arbitrary  connected Lie group and $\mathfrak{g} $ its Lie algebra. Let $\{\xi_{1},\ldots,\xi_{l}\}$ and
$\left\{  \epsilon^{1},\ldots,\epsilon^{l}\right\}  $ be dual bases of
$\mathfrak{g}$ and $\mathfrak{g}^{\ast}$, respectively.  Left (respectively, right) translations on $G$ will be denoted by $L:G\times G\rightarrow G$ (respectively, $R:G\times G\rightarrow G$). 
With the same notation that we have used so far, let
\begin{equation}%
\begin{array}
[c]{rrl}%
S\left(  \mu,g\right)  :T_{\mu}\mathfrak{g}\simeq\mathfrak{g} &
\longrightarrow & T_{g}G\\
\eta & \longmapsto & \sum_{i=1}^{l}\xi_{i}^{G}\left(  g\right)  \left\langle
\epsilon^{i},\eta\right\rangle =\eta^{G}(g)
\end{array}
\label{eq 8 bis}%
\end{equation}
be a Stratonovich operator from $\mathfrak{g}$ to $G$, where $\eta^{G}$
denotes the infinitesimal generator associated to the $G$-action on itself by
left translations. Consider the stochastic differential equation associated
to (\ref{eq 8 bis}),%
\begin{equation}
\delta g_{t}=\sum_{i=1}^{l}\xi_{i}^{G}\left(  g_t\right)  \delta X_{t}^{i},
\label{eq 17}%
\end{equation}
for some driving noise (semimartingale) $X:\mathbb{R}_{+}\times\Omega\rightarrow\mathfrak{g}$.
Using the equivariance of the vector fields $\xi^{G}\in\mathfrak{X}(G)$ with
respect to right translations, that is, $T_{h}
R_{g}(\xi^{G}(h)))=\xi^{G}(R_{g}(h))$ for any $g$, $h\in G$, and $\xi
\in\mathfrak{g}$, it is immediate to check that if $\Gamma^{e}$ is the
solution of (\ref{eq 17}) with initial condition $\Gamma_{t=0}^{e}=e$ a.s.,
then the solution $\Gamma_t^g$ starting at $g\in G$ is given by
\begin{equation}
\Gamma_{t}^{g}=L_{\Gamma_{t}^{e}}g=R_{g}\left(  \Gamma_{t}^{e}\right)
\label{eq 28}%
\end{equation}
In other words, the stochastic differential equation~(\ref{eq 17}) {\it has a superposition rule in the sense of Definition~\ref{def 1} and the superposition function $\Phi$ is
given by}
\[%
\begin{array}
[c]{rrl}%
\Phi:G\times G & \longrightarrow & G\\
\left(  h,g\right)  & \longmapsto & L_{h}g=R_{g}h.
\end{array}
\]
It is also worth noticing that (\ref{eq 17}) is {\bfi stochastically complete}
(\cite[Chapter VII \S 6]{elworthy book}) since it is a left-invariant system.
Therefore any solution of (\ref{eq 17}) is defined for all $\left(
t,\omega\right)  \in\mathbb{R}_{+}\times\Omega$ and, consequently, so is any
one point motion and, in particular, any solution of any Lie-Scheffers system on a manifold
$M$ which can be globally considered as induced by a group action $\Xi:G\times
M\rightarrow M$.

\medskip

\noindent {\bf L\'evy processes and Lie-Scheffers systems.} This is an important class of Lie group valued stochastic processes and, as we will now see, a class of examples of Lie-Scheffers systems. Recall that a continuous
process $g:\mathbb{R}_{+}\times\Omega\rightarrow G$ is called a right L\'{e}vy
process if, for any $0=t_{0}<t_{1}<t_{2}<\ldots<t_{n}$, the increments%
\begin{equation}
g_{t_{0}},g_{t_{0}}g_{t_{1}}^{-1},g_{t_{1}}g_{t_{2}}^{-1},\ldots,g_{t_{n-1}%
}g_{t_{n}}^{-1} \label{eq 42}%
\end{equation}
are \emph{independent }and\emph{\ stationary}. This means that the 
random variables in (\ref{eq 42}) are mutually independent and that their distributions only depend on the differences $t _i-t_{i-1} $, $ i \in \{ 1, \ldots, n\} $. If $g_{t_{0}}\neq e$ a.s., we define $g_{t}^{e}%
=g_{t}g_{t_{0}}^{-1}$, which is a right L\'{e}vy process starting at the
identity. 

We are now going to see that continuous L\'{e}vy processes and Lie-Scheffers systems are closely
related. First of all, recall that any right
L\'{e}vy process on a locally compact topological group with a countable basis
of open sets is a Markov process with a right invariant Feller transition
semigroup $\left\{  P_{t}\right\}  _{t\in\mathbb{R}_{+}}$ given by
$P_{t}f\left(  g\right)  :=E\left[  f\left(  g_{t}^{e}g\right)  \right]  $,
$g\in G$, where $f:G\rightarrow\mathbb{R}$ is any measurable function.
Conversely, any right invariant \emph{continuous} Markov process is a right
L\'{e}vy process (\cite[Proposition 1.2]{Liao book}). Moreover, if
$g:\mathbb{R}_{+}\times\Omega\rightarrow G$ is a right L\'{e}vy process, then there
exists a $l$-dimensional Brownian motion $B:\mathbb{R}_{+}\times
\Omega\rightarrow\mathbb{R}^{l}$ with respect to the natural filtration
$\{\mathcal{F}_{t}^{e}\}_{t\in\mathbb{R}_{+}}$ of the process $g_{t}^{e}$,
$l=\dim\left(  \mathfrak{g}\right)  $, with covariance matrix $(a_{ij}%
)_{i,j=1,\ldots,l}$ and some constants $\{c_{i}\}_{i=1,\ldots,l}$ such that%
\[
f\left(  g_{t}\right)  =f\left(  g_{0}\right)  +\sum_{i=1}^{l}\int_{0}^{t}%
\xi_{i}^{G}[f]\left(  g_{s}\right)  \delta B_{s}^{i}+\sum_{i=1}^{l}c_{i}%
\int_{0}^{t}\xi_{i}^{G}[f]\left(  g_{s}\right)  ds,
\]
for any $f\in C^{2}\left(  G\right)  $ and where, as before, $\left\{  \xi
_{1},\ldots,\xi_{l}\right\}  $ is a basis of $\mathfrak{g}$ (\cite[Theorem
1.2]{Liao book}). This expression amounts to saying that the L\'evy process $g:\mathbb{R}_{+}\times\Omega\rightarrow G$ satisfies the stochastic differential equation
\[
\delta g_{t}  = \sum_{i=1}^{l}c_{i}
\xi_{i}^{G}\left(  g_{s}\right)  \delta s+\sum_{i=1}^{l}
\xi_{i}^{G}\left(  g_{s}\right)  \delta B_{s}^{i},
\]
and hence by Corollary~\ref{corolario Lie clasico} we can conclude that
\emph{any continuous right L\'{e}vy process is a
solution of a right invariant Lie-Scheffers system}. Additionally, it can be shown in this context (see~\cite[Theorem 1.2]{Liao book}) that one point motions obtained out of a $G$-action
$\Xi:G\times M\rightarrow M$ are Markov processes with Feller transition
semigroup $\left\{  P_{t}^{M}\right\}  _{t\in\mathbb{R}_{+}}$%
\[
P_{t}^{M}f\left(  z\right)  =E\left[  f(\Xi\left(  g_{t}^{e},z\right)
)\right]  ,~z\in M,~f\in C\left(  M\right).  
\]

\medskip

\noindent {\bf Lie-Scheffers systems on homogeneous spaces.} Let $H\subset G$ be a closed subgroup of $G$ and consider the
{\bfi homogeneous space} $G/H=\left\{  gH\mid g\in G\right\}  $ with the
unique smooth structure that makes the projection $\pi_{H}:G\rightarrow G/H$
into a submersion. The group $G$ acts on $G/H$ via the map $\lambda:G\times
G/H\rightarrow G/H$ on $G/H$ defined by $\left(  h,gH\right)  \mapsto(hg)H$.
It is immediate to check that the infinitesimal generators associated to the left
$G$-actions on $G$ and on $G/H$ are $\pi_{H}$-related, that is,
\[
T_{g}\pi_{H}\left(  \xi^{G}\left(  g\right)  \right)  =\xi^{G/H}\left(
\pi_{H}\left(  g\right)  \right)
\]
for any $g\in G$, any $\xi\in\mathfrak{g}$, and where $\xi^{G/H}\left(
gH\right)  =\left.  \frac{d}{dt}\right\vert _{t=0}\lambda_{\exp\left(
t\xi\right)  }\left(  gH\right)  $. This straightforward observation has as an
immediate consequence the next proposition:

\begin{proposition}
\label{proposicion reduccion}Let $X:R_{+}\times\Omega\rightarrow\mathfrak{g}$
be a $\mathfrak{g}$-valued semimartingale, $G$ a Lie group, and $H\subset G$ a
closed subgroup. Let $\Gamma$ be a solution of the Lie-Scheffers system
defined by $X$ and the Stratonovich operator (\ref{eq 8 bis}) with initial
condition $\Gamma_{t=0}$. Then, $\pi_{H}\left(  \Gamma\right)  $ is a solution
of the Lie-Scheffers system on $G/H$
\begin{equation}
\delta\overline{\Gamma}=\sum_{j=1}^{l}\xi_{j}^{G/H}\left(  \overline{\Gamma
}_{t}\right)  \delta X_{t}^{j} \label{eq 10}%
\end{equation}
with initial condition $\pi_{H}\left(  \Gamma_{t=0}\right)  $.
\end{proposition}

Observe that since the Stratonovich operator~(\ref{eq 8 bis}) is right
invariant by the action of $G$, and therefore $H$-invariant, and that since
this action is free and proper, the previous proposition can be seen as a
particular case of the {\bfi Reduction Theorem} in \cite{paper reduccion}. The
next theorem is a transcription of the {\bfi  Reconstruction Theorem}
in~\cite{paper reduccion} into the present context and describes how to
construct solutions in the opposite direction, that is, it tells us \emph{how
to construct a solution $\Gamma$ of the Lie-Scheffers system (\ref{eq 17}) out
of the solutions of two other dimensionally smaller Lie-Scheffers systems}:
first, a solution of the reduced system~(\ref{eq 10}) and second, another
solution of a new Lie-Scheffers system, now on $H$.

\begin{theorem}
Let $X:\mathbb{R}_{+}\times\Omega\rightarrow\mathfrak{g}$ be a $\mathfrak{g}
$-valued semimartingale, $G$ a Lie group, $H\subset G$ a closed subgroup, and
$S$ the Stratonovich operator defined in (\ref{eq 8 bis}). Let $R:H\times
G\rightarrow G$ be the (right) action of $H$ on $G$ by right translations and
$A$ an auxiliary principal connection on $\pi_{H}:G\rightarrow G/H$. Then, any
solution $\Gamma$ of the system (\ref{eq 17}) can be written in the form%
\[
\Gamma_{t}=R_{h_{t}}g_{t}=g_{t}h_{t}.
\]
In this statement, $g:\mathbb{R}_{+}\times\Omega\rightarrow G$ is a $G$-valued
semimartingale horizontal with respect to $A$, i.e. $\int\left\langle A,\delta
g_{t}\right\rangle =0\in\mathfrak{g}$, $g_{t=0}=\Gamma_{t=0}$, and such that
$\pi_{H}\left(  g_{t}\right)  $ is a solution of the reduced system
(\ref{eq 10}). On the other hand, $h:\mathbb{R}_{+}\times\Omega\rightarrow H$
is a $H$-valued semimartingale that satisfies the stochastic differential
equation%
\begin{equation}
\delta h_{t}=\widetilde{R}\left(  Y_{t},h_{t}\right)  \delta Y_{t}
\label{eq 25 bis}%
\end{equation}
with initial condition $h_{t=0}=e$, and associated to the Stratonovich
operator
\begin{equation}%
\begin{array}
[c]{rrl}%
\widetilde{R}(\xi,h):T_{\xi}\mathfrak{h} & \longrightarrow & T_{h}H\\
\eta & \longmapsto & T_{e}R_{h}(\eta)= \eta^{H}(h),
\end{array}
\label{eq 25}%
\end{equation}
and the stochastic component $Y:\mathbb{R}_{+}\times\Omega\rightarrow
\mathfrak{h}$ given by%
\[
Y=\sum_{i=1}^{l}\int A_{g_{t}}\left(  \xi_{i}^{G}\left(  g_{t}\right)
\right)  \delta X^{i}.
\]

\end{theorem}

\noindent\textbf{Proof.\ \ } See \cite[Theorem 3.2 and Proposition 3.4]{paper
reduccion}. \quad$\blacksquare$

\subsection{The Wei-Norman method for solving stochastic Lie-Scheffers
systems}

\label{wei norman method}

The method that we are going to develop in this subsection is a generalization
to stochastic systems of the one proposed by Wei and Norman in
\cite{wei-norman1, wei-norman2} in order to solve by quadratures time
evolution equations of the form $\frac{dU_{t}}{dt}=H_{t}U_{t}$ that appear in
quantum mechanics, where both $U_{t}$ and $H_{t}$ are bounded linear operators
on a suitable Hilbert space. This method has already been adapted by
Cari\~{n}ena and Ramos~\cite{arturo ramos} to the study of deterministic
Lie-Scheffers systems on Lie groups and it is their approach that we will
follow. As we will see later on, the power of this method and the ease of its
implementation depends strongly on the algebraic structure of the Lie algebra
$\mathfrak{g}$ of the group $G$ where the solutions of the stochastic
differential equation take values.

Let $\Gamma:\mathbb{R}\times\Omega\rightarrow G$ be the solution of
(\ref{eq 17}) such that $\Gamma_{t=0}=e\in G$ a.s.; we write it down in terms
of second kind canonical coordinates with respect to a basis $\{\xi_{1}%
,\ldots,\xi_{l}\}$ of the Lie algebra $\mathfrak{g}$. That is,%
\begin{equation}
\Gamma_{t}=\exp(d_{t}^{1}\xi_{1})\cdots\exp(d_{t}^{l}\xi_{l}), \label{eq 32}%
\end{equation}
where $\{d_{t}^{1},\ldots,d_{t}^{l}\}$ is a family of real-valued
semimartingales, $d^{i}:\mathbb{R}_{+}\times\Omega\rightarrow\mathbb{R}$, such
that $d_{t=0}^{i}=0$ a.s. for any $i=1,\ldots,l$. Notice that the expression
(\ref{eq 32}) is only valid up to the exit time of $\Gamma$ from the
neighborhood $U_{e}$ of $e\in G$ where the second kind canonical coordinates
for $G$ around the origin are valid. The key idea in this method is that if
the functions $d^{i}$ were differentiable then%
\[
\frac{d\Gamma_{t}}{dt}=T_{e}R_{\Gamma_{t}}\left(  \sum\nolimits_{i=1}^{l}%
\dot{d}_{t}^{i}\left(
{\displaystyle\prod\nolimits_{j<i}}
\operatorname*{Ad}\nolimits_{\exp\left(  d_{t}^{j}\xi_{j}\right)  }\right)
\xi_{i}\right)
\]
(see \cite[Eq. (33) and (34)]{arturo ramos}), where $\operatorname*{Ad}%
\nolimits_{g}(\eta)\in\mathfrak{g}$ is the adjoint representation of $G$ on
$\mathfrak{g}$, $g\in G$, $\eta\in\mathfrak{g}$. In our setup we obviously cannot  invoke the differentiability of the functions $d ^i$, however  applying the
Stratonovich differentiation rules to (\ref{eq 32}) with $d^{i}$ our
real-valued semimartingales, $i=1,\ldots,l$, we have%
\[
\delta\Gamma_{t}=T_{e}R_{\Gamma_{t}}\left(  \sum\nolimits_{i=1}^{l}\delta
d_{t}^{i}\left(
{\displaystyle\prod\nolimits_{j<i}}
\operatorname*{Ad}\nolimits_{\exp\left(  d_{t}^{j}\xi_{j}\right)  }\right)
\xi_{i}\right)  .
\]
This expression implies that for any right invariant one-form $\mu^{G}%
\in\Omega(G)$, that is, $\mu^{G}(g)=T_{g}^{\ast}R_{g^{-1}}(\mu)$ for any $g\in
G$ and a fixed $\mu\in\mathfrak{g}^{\ast}$,
\begin{equation}
\int\left\langle \mu^{G},\delta\Gamma\right\rangle =\langle\mu,\sum_{i=1}%
^{r}\int\left(
{\displaystyle\prod\nolimits_{j<i}}
\operatorname*{Ad}\nolimits_{\exp\left(  \sum_{j=1}^{l}d_{t}^{j}\nu
_{j}\right)  }\right)  \xi_{i}\delta d_{t}^{i}\rangle.
\label{pair with g invariant}%
\end{equation}
At the same time, it is clear that $\int\left\langle \mu^{G},\delta
\Gamma\right\rangle =\left\langle \mu,X\right\rangle $ and
hence~(\ref{pair with g invariant}) implies that
\[
X=\sum_{i=1}^{l}\int\left(
{\displaystyle\prod\nolimits_{j<i}}
\operatorname*{Ad}\nolimits_{\exp\left(  d_{t}^{j}\xi_{j}\right)  }\right)
\xi_{i}\delta d_{t}^{i}.
\]
Using the identity $\operatorname*{Ad}\nolimits_{\exp(\eta)}=\operatorname*{e}%
^{\operatorname*{ad}(\eta)}=\sum_{n\geq0}\frac{1}{n!}\operatorname*{ad}%
(\eta)\circ\overset{n}{\ldots}\circ\operatorname*{ad}(\eta)$, for any $\eta
\in\mathfrak{g}$, and writing $X=\sum_{i=1}^{l}X^{i}\xi_{i}$, we get the
relation
\begin{equation}
\sum_{i=1}^{l}X^{i}\xi_{i}=\sum_{i=1}^{l}\int\left(
{\displaystyle\prod\nolimits_{j<i}}
\operatorname*{e}{}^{\operatorname*{ad}\left(  d_{t}^{j}\xi_{j}\right)
}\right)  \xi_{i}\delta d_{t}^{i}. \label{eq 33}%
\end{equation}
The system of stochastic differential equations (\ref{eq 33}) can be solved
for the semimartingales $d_{t}^{i}$, $i=1,\ldots,m$ by quadratures if the Lie
algebra $\mathfrak{g}$\ is solvable (see \cite{wei-norman1, wei-norman2}) and,
in particular, for nilpotent Lie algebras. The solvable case was extensively
studied in \cite{Kunita seminaire} where similar conclusions were presented
using a different approach.

As a simple example consider the affine group in one dimension $\mathcal{A}%
_{1}$,that is, the group of affine transformations of the real line. Any
element of $\mathcal{A}_{1}$ can be expressed as a pair of real numbers
$\left(  a_{0},a_{1}\right)  $ with $a_{1}\neq0$ defining the affine
transformation $x\mapsto a_{1}x+a_{0}$. The product $\ast:\mathcal{A}%
_{1}\times\mathcal{A}_{1}\rightarrow\mathcal{A}_{1}$ in $\mathcal{A}_{1} $ is
\[
\left(  a_{0},a_{1}\right)  \ast\left(  b_{0},b_{1}\right)  =\left(
a_{0}+a_{1}b_{0},a_{1}b_{1}\right)  .
\]
If $\left\{  \xi_{0}=\left(  1,0\right)  ,\xi_{1}=\left(  0,1\right)
\right\}  $ is a basis of the Lie algebra $\mathfrak{a}_{1}$ of $\mathcal{A}%
_{1} $, it is immediate to check that
\begin{equation}
\left[  \xi_{0},\xi_{1}\right]  =\operatorname*{ad}\nolimits_{\xi_{0}}(\xi
_{1})=-\xi_{0}. \label{eq 39}%
\end{equation}
Furthermore, the infinitesimal generators associated to the left action of
$\mathcal{A}_{1}$ on itself are%
\[
\xi_{0}^{\mathcal{A}_{1}}\left(  x,y\right)  =\frac{\partial}{\partial
x}\text{ \ \ and \ \ }\xi_{1}^{\mathcal{A}_{1}}\left(  x,y\right)
=x\frac{\partial}{\partial x}+y\frac{\partial}{\partial y}.
\]
A typical Lie system on $\mathcal{A}_{1}$ would be, for instance, the
following Stratonovich differential equation on the upper half-plane
$H_{+}=\left\{  \left(  x,y\right)  \in\mathbb{R}^{2}~|~y>0\right\}  $,%
\[
\delta\Gamma_{x}=dt+\Gamma_{x}\delta B_{t},~~\delta\Gamma_{y}=dt+\Gamma
_{y}\delta B_{t}%
\]
obtained as a particular case of~(\ref{eq 8 bis}) when $G= \mathcal{A}_{1} $
and $X=(t, B)$, where $B:\mathbb{R}_{+}\times\Omega\rightarrow\mathbb{R}$ is a
Brownian motion. More generally, let $X:\mathbb{R}_{+}\times\Omega
\rightarrow\mathfrak{a}_{1}$ be an $\mathfrak{a}_{1}$-valued semimartingale
and write $X=X^{0} \xi_{0}+X^{1}\xi_{1}$, with $X^{0}$ and $X^{1}$ real
semimartingales. Then, using (\ref{eq 39}), (\ref{eq 33}) reads in this
particular case
\[
X^{0}\xi_{0}+X^{1}\xi_{1}=\int\xi_{0} \delta d_{t} ^{0}+\int\left(  \xi_{1}-
d_{t} ^{0} \xi_{0}\right)  \delta d^{1}_{t}= \left(  \int\delta d_{t}^{0}-\int
d_{t} ^{0}\delta d_{t}^{1}\right)  \xi_{0} +\left(  \int\delta d_{t}%
^{1}\right)  \xi_{1}.
\]
Putting together the terms that go both with $\xi_{1}$ and $\xi_{0}$
respectively, we obtain
\[
d_{t}^{1} =X_{t}^{1}, \quad d_{t}^{0} =X_{t}^{0}+\int_{0}^{t}d_{s}^{0}\delta
X_{s}^{1},
\]
and hence
\[
d_{t}^{0}=\operatorname*{e}\nolimits^{X_{t}^{1}}\left(  \int_{0}^{t}\delta
X_{s}^{0}\operatorname*{e}\nolimits^{-X_{s}^{1}}\right)  .
\]

\section{The flow of a stochastic Lie-Scheffers system}
\label{The flow of a stochastic Lie-Scheffers system}

Theorem \ref{teorema Lie} claims, roughly speaking,  that the stochastic system
(\ref{stochastic differential equation expression}) admits a superposition
rule $\left(  \Phi,\{\Gamma_{1},\ldots,\Gamma_{m}\}\right)  $ if the
components of the Stratonovich operator $S\left(  x,z\right)
:T_{x}\mathbb{R}^{l}\longrightarrow T_{p}\mathbb{R}^{n}$, $x\in\mathbb{R}^{l}
$, $p\in\mathbb{R}^{n}$, that define it may be written as $S_{j}(X,z)=\sum_{i=1}^{r}b_{j}%
^{i}\left(  X\right)  Y_{i}\left(  z\right)  $, where $b_{j}^{i}\in C^{\infty
}(\mathbb{R}^{l})$ and $\{Y_{1},\ldots,Y_{r}\}\subseteq\mathfrak{X}\left(
\mathbb{R}^{n}\right)  $ span an involutive distribution. The converse of this statement is also true
provided that, for a given initial condition
$z\in\mathbb{R}^{n}$, the point $\left(  z,(\Gamma_{1},\ldots,\Gamma
_{m})_{t=0}\right)  $ is a regular point of the foliation $\mathcal{G}_{0}$
generated by the diagonal extensions of $\{S_{1}(X,\cdot),\ldots,S_{m}%
(X,\cdot)\}$. Notice that this is a reasonable condition
since the set of regular points of a  generalized foliation is open and dense (\cite[Th\'eor\`eme 2.2]{dazord 1985}). Moreover, when this happens, the vector fields $\{Y_{1},\ldots
,Y_{r}\}$ form a real Lie algebra. 

The condition on the vector fields  $\{Y_{1},\ldots
,Y_{r}\}$ forming a real finite dimensional Lie algebra or, more generally, $\dim\left(
\mathrm{Lie}\{Y_{1},\ldots,Y_{r}\}\right)  <\infty$, are particularly appealing since these are algebraic requirements that we may expect to be easily verified for stochastic
differential equations of a certain type. 
Moreover, these conditions have consequences that go beyond
Corollary~\ref{corolario Lie clasico}.
More especifically, we will show that if $\dim\left(  \mathrm{Lie}\{Y_{1},\ldots,Y_{r}\}\right)
<\infty$, then the general solution of a stochastic differential equation {\it can be
written by composing a deterministic function with a suitable noise}. In the following paragraphs we are going to give a precise meaning to this statement and to put it in the context of well known results available in the literature. 

Traditionally, stochastic differential equations on a manifold $M$ have been
presented as
\begin{equation}
\delta\Gamma_{t}=Y_{0}(\Gamma_{t})dt+\sum_{i=1}^{r}Y_{i}\left(  \Gamma
_{t}\right)  \delta B_{t}^{i}, \label{eq 37}%
\end{equation}
where $\left\{  Y_{0},\ldots,Y_{r}\right\}  \subseteq\mathfrak{X}\left(
M\right)  $ and $B:\mathbb{R}_{+}\times\Omega\rightarrow\mathbb{R}^{r}$ is a
$r$-dimensional Brownian motion defined on a standard filtered probability
space $(\Omega, \mathcal{F}_t, P)$. For the sake of having a more compact notation, we write$B_{t}^{0}:=t$. The flow of such a stochastic differential equation may be locally written, that is, up to a given stopping time $\tau$, by means of a Taylor
series expansion that comes out of Picard's iterative method for solving stochastic differential
equations. In order to be more explicit we
introduce some notation. Let $J=\{j_{1},\ldots,j_{n}\}$, $j_{i}
\in\{0,\ldots,r\}$, $1\leq i\leq n$, be a multi-index of {\bfi size} $n$.
$\left\Vert J\right\Vert $ will denote the {\bfi degree} of $J$ that, by definition, is the
size of $J$ plus the number of zeros in the $n$-tuple $(j_{1},\ldots,j_{n}) $. For any $J=\{j_{1},\ldots,j_{n}\}$, we consider the
iterated Stratonovich multiple integral
\[
B_{t}^{J}=\idotsint\limits_{0<t_{1}<\ldots<t_{n}<t}\delta B_{t_{1}}^{j_{1}%
}\cdots\delta B_{t_{n}}^{j_{n}}.
\]
In addition, $Y_{J}$ will denote
\[
Y_{J}:=[Y_{j_{1}},[Y_{j_{2}},\ldots,\lbrack Y_{j_{n-1}},Y_{j_{n}}]].
\]
If $Y\in\mathfrak{X}\left(  M\right)  $ is a vector field on the
manifold $M$, we will use the following notation for its flow: $\exp\left(  sY\right)  (z)$ denotes the solution at time $s$ of
the ordinary differential equation $\dot{\gamma}=Y(\gamma)$ with initial
condition $\gamma(0)=z$. Then,

\begin{theorem}
[{\cite[Th\'{e}or\`{e}me 20]{Ben Arous}}]\label{teorema Ben Arous}With the
notation introduced so far, if $\dim\left(  \operatorname*{Lie}\{Y_{0}%
,\ldots,Y_{r}\}\right)  <\infty$ and $\operatorname{span}\{\operatorname*{Lie}
\{Y_{0},\ldots,Y_{r}\}\}$ has constant dimension on a neighborhood
$V$ of the point $z\in  M$, then there exists a stopping time $\tau$ such that the solution of
(\ref{eq 37}) with initial condition $z$ can be expressed as%
\begin{equation}
\Gamma_{t}^{z}=\exp\left(  \sum_{n=1}^{\infty}\sum_{\left\Vert J\right\Vert
=n}\beta_{J}B_{t}^{J}\right)  (z) \label{eq 43}%
\end{equation}
up to time $\tau$. In this expression, $\beta_{J}:=\sum_{\sigma\in S_{n}}%
\frac{(-1)^{e(\sigma)}}{n^{2}\binom{n-1}{e(\sigma)}}Y_{\sigma(J)}$, $S_{n}$
denotes the permutation group of $n$ elements, and $e(\sigma)$ is the
cardinality of the set $\{j\in\{1,\ldots,m-1\}~|~\sigma(j)>\sigma(j+1)\}$.
\end{theorem}

If the finiteness condition on the dimensionality of the Lie algebra generated by the vector fields is not available but, nevertheless, $\{Y_{0},\ldots,Y_{r}\}$ are Lipschitz vector
fields, then the solution of (\ref{eq 37}) starting at $z\in M$ can always be
approximated by a process like (\ref{eq 43}): if $\zeta_{t}^{N}$ denotes the
\emph{finite} sum $\sum_{n=1}^{N}\sum_{\left\Vert J\right\Vert =n}\beta
_{J}B_{t}^{J}$, then%
\[
\Gamma_{t}^{z}=\exp\left(  \zeta_{t}^{N}\right)  (z)+t^{N/2}R_{N}(t)
\]
where the error term $R_{N}(t)$ is bounded in probability when $t$ tends to
$0$ (\cite[Theorem 2.1]{Castell - asymptotic expansion}). The expression (\ref{eq 43}) also
holds if instead of the hypotheses of Theorem \ref{teorema Ben Arous} we
require $M$ to be an analytic manifold and $\{Y_{0},\ldots,Y_{r}\}$ a family
of real analytic vector fields (\cite[Th\'{e}or\`{e}me 10]{Ben Arous}).
An important consequence of Theorem \ref{teorema Ben Arous} lies in the
fact that the general solution of the stochastic differential equation
(\ref{eq 37}) may be written, at least locally and up to a suitable stopping
time $\tau$, as the composition of a deterministic and smooth function, namely,
the flow exponential, with the diffusion that defines the stochastic differential equation (see \cite{Hu seminaire}
for a complementary reading). 
From this point of view, there is a strong resemblance between Theorem~\ref{teorema Ben Arous} and Theorem~\ref{teorema Lie}: 
\begin{itemize}
\item First, by Corollary \ref{corolario Lie clasico}, all the systems that satisfy the hypotheses of Theorem \ref{teorema Ben Arous} admit a superposition rule. 
\item Second, the superposition rule allows us to write any solution as a composition of the deterministic function and the set of solutions $\left\{  \Gamma_{1},\ldots,\Gamma_{m}\right\}  $ that are responsible for the stochastic behavior of the resulting flow.
\end{itemize}

We conclude by quoting two references that study the nilpotent
case (that is, the Lie algebra $\mathrm{Lie}\{Y_{0},\ldots,Y_{r}\}$ is
nilpotent); this case has deserved special attention in the literature (see, for example, \cite{Kunita seminaire}) because in that situation
the Taylor series expansion of the
flow in terms of iterated integrals in (\ref{eq 43}) becomes finite. We also
recommend the excellent exposition in \cite{baudoin} for a complementary approach to the subject of Taylor series approximation of the general solution of (\ref{eq 37}); in this book it is shown that, for instance,  the Carnot
group of depth $N=\dim\left(  \mathrm{Lie}\{Y_{0},\ldots,Y_{r}\}\right)  $ can
be used in the nilpotent case to integrate the Lie algebra action of
$\mathrm{Lie}\{Y_{0},\ldots,Y_{r}\}$ when one writes, as we did in the previous section, a Lie-Scheffers system as a stochastic differential equation on a Lie group that acts on the manifold in question.

\section{Examples.}

\subsection{Inhomogeneous linear systems.}

Let $A_{k}:\mathbb{R}\rightarrow M_{n}(\mathbb{R})$ a $n\times n$
time-dependent real matrix and $B_{k}:\mathbb{R}\rightarrow\mathbb{R}^{n}$ a
time-dependent vector for any $k=1,\ldots,l$. Let $X:\mathbb{R}_{+}%
\times\Omega\rightarrow\mathbb{R}^{l}$ be a semimartingale. An inhomogeneuous
linear system is a system of stochastic differential equations on
$\mathbb{R}^{n}$ that may be written as%
\begin{equation}
\delta\Gamma_{t}=\sum_{k=1}^{l}(A_{k}(t)(\Gamma_{t})-B_{k}(t))\delta X_{t}^{k}
\label{eq 34}%
\end{equation}
Let $\left(  q^{1},\ldots,q^{n}\right)  $ be coordinates for $\mathbb{R}^{n}$.
It is an exercise to check that (\ref{eq 34}) can be equivalently written as%
\[
\delta\Gamma_{t}=\sum_{k=1}^{l}\sum_{i,j=1}^{n}\left(  A_{k}\right)  _{i}%
^{j}(t)Y_{j}^{i}(\Gamma_{t})\delta X_{t}^{k}+\sum_{k=1}^{l}\sum_{i,j=1}%
^{n}(B_{k})^{j}(t)Z_{j}(\Gamma_{t})\delta X_{t}^{k}%
\]
where the vector fields $Y_{j}^{i}$, $Z_{j}\in\mathfrak{X}(\mathbb{R}^{n})$,
$i,j,k=1,\ldots,n$, are given by
\[
Y_{j}^{i}=q^{i}\frac{\partial}{\partial q^{j}}\text{, \ \ \ }Z_{j}%
=\frac{\partial}{\partial q^{j}}.
\]
Given that
\[
\lbrack Y_{j}^{i},Y_{l}^{k}]=\delta_{j}^{k}Y_{l}^{i}-\delta_{l}^{i}Y_{j}%
^{k}\text{,\ \ \ }[Y_{j}^{i},Z_{k}]=-\delta_{k}^{i}Z_{j}\text{, \ and
\ }[Z_{i},Z_{j}]=0
\]
we see that the vectors $\{Y_{j}^{i},Z_{k}~|~i,j,k=1,\ldots,n\}\subset
\mathfrak{X}(\mathbb{R}^{n})$ span a Lie algebra isomorphic to the $(n^{2}%
+n)$-dimensional Lie algebra of the group of affine transformations of
$\mathbb{R}^{n}$. Therefore, the system (\ref{eq 34}) satisfies the hypotheses
of Theorem \ref{teorema Lie} and hence it admits a superposition rule. In
order to explicitly construct the superposition rule, let $\Gamma^{e_{j}}$ be
the solution of the homogeneous part of (\ref{eq 34}),%
\[
\delta\Gamma_{t}=\sum_{k=1}^{l}A_{k}(t)(\Gamma_{t})\delta X_{t}^{k}%
\]
with initial solution $\Gamma_{t=0}^{e_{j}}=e_{j}\in\mathbb{R}^{n}$ a.s.,
where $e_{j}=(0,\overset{j-1}{\ldots},0,1,0,\ldots,0)$ for any $j=1,\ldots,n$.
Let $\overline{\Gamma}$ be a particular solution of (\ref{eq 34}) with initial
condition $\overline{\Gamma}_{t=0}=0\in\mathbb{R}^{n}$ a.s.. Then,%
\[
\Gamma_{t}=\sum_{j=1}^{n}z^{j}\Gamma_{t}^{e_{j}}+\overline{\Gamma}_{t}%
\]
is the general semimartingale solution of (\ref{eq 34}) starting at
$z=(z^{1},\ldots,z^{n})\in\mathbb{R}^{n}$.

\subsection{The stochastic exponential of a Lie group.}

Let $G$ be a Lie group and $\mathfrak{g}$ its Lie algebra. Let $\{\xi
_{1},\ldots,\xi_{l}\}$ a basis of $\mathfrak{g}$ and $X:\mathbb{R}_{+}%
\times\Omega\rightarrow\mathfrak{g}$ be a $\mathfrak{g}$-valued
semimartingale. Observe that $X$ can be written as $X=\sum_{i=1}^{r}a_{t}%
^{i}\xi_{i}$ for a family of real semimartingales $a^{i}:\mathbb{R}_{+}%
\times\Omega\rightarrow\mathbb{R}$, $i=1,\ldots,l$. Following \cite{hakim} and
\cite{BSDE Lie groups}, we define the (left) {\bfi stochastic exponential}
$\mathcal{E}(X):\mathbb{R}_{+}\times\Omega\rightarrow G$ of $X$ as the unique
solution of the Lie-Scheffers system on $G$ given by%
\[
\delta\Gamma_{t}=\sum_{i=1}^{l}(\xi_{i})^{G}(\Gamma_{t})\delta a_{t}^{i}%
\]
with initial condition $\Gamma_{t=0}=e\in G$ a.s.. Unlike the conventions used
in Section~\ref{Lie-Scheffers systems on Lie groups and homogeneous spaces},
the vector fields $(\xi_{i})^{G}\in\mathfrak{X}(G)$ here are not the
right-invariant vector fields built from $\xi_{i}$, $i=1,\ldots,l$, but the
left-invariant ones. That is,%
\[
(\xi_{i})^{G}(g)=T_{e}L_{g}(\xi_{i}),~~g\in G.
\]
Except for the fact that $(\xi_{i})^{G}\in\mathfrak{X}(G)$, $i=1,\ldots,l$,
are now left-invariant, solving a Lie-Scheffers system on a Lie group such as
those presented in
Section~\ref{Lie-Scheffers systems on Lie groups and homogeneous spaces}
amounts to computing the stochastic exponential of a given $\mathfrak{g}%
$-valued semimartingale $X$.

The stochastic exponential establishes a bijection between $\mathfrak{g}%
$-valued local martingales and martingales on $G$ with respect to certain
connections. Recall that, given an affine connection $\nabla:\mathfrak{X}%
(M)\times\mathfrak{X}(M)\rightarrow\mathfrak{X}(M)$ on a manifold $M$, a
$M$-valued semimartingale $\Gamma:\mathbb{R}_{+}\times\Omega\rightarrow M$ is
said to be a $\nabla$-martingale (or a martingale with respect to $\nabla$)
provided that
\[
f(\Gamma)-f(\Gamma_{t=0})-\frac{1}{2}\int\operatorname*{Hess}f\left(
d\Gamma,d\Gamma\right)
\]
is a real local martingale for any $f\in C^{\infty}(M)$, where
$\operatorname*{Hess}f:\mathfrak{X}(M)\times\mathfrak{X}(M)\rightarrow
C^{\infty}(M)$ is the bilinear form defined as%
\[
\operatorname*{Hess}f\left(  Y,Z\right)  =Y\left[  Z\left[  f\right]  \right]
-\nabla_{Z}Y\left[  f\right]
\]
for any $Y$, $Z\in\mathfrak{X}(M)$ (see \cite[Chapter IV]{emery}). When $M=G$
is a Lie group, one can construct left invariant connections $\nabla$ by using
bilinear skew-symmetric forms $\alpha:\mathfrak{g\times g\rightarrow
}\mathbb{R}$ on the Lie algebra $\mathfrak{g}$ via the definition
\[
\nabla_{ \xi^{G} }\eta^{G}:=\alpha(\xi,\eta),~~\xi,\eta\in\mathfrak{g.}%
\]
The curves $\exp(t\xi)\in G$, where $\xi\in\mathfrak{g}$ and $\exp
:\mathfrak{g}\rightarrow G$ is the Lie algebraic exponential, coincide with
the geodesics $c (t) $ with respect to these connections that start at $e\in
G$ and that satisfy $\dot c (0)= \xi$. It can be shown (\cite[Lemma 1.4]{BSDE
Lie groups}) that the connections built from $\alpha=0$ and $\alpha(\xi
,\eta)=\frac{1}{2}\left[  \xi,\eta\right]  $ induce the same $\nabla
$-martingales on $G$. Moreover, with respect to these two connections, the set
of $\nabla$-martingales consists precisely of the processes of the form
$\Gamma_{0}\mathcal{E}(X)$ where $X$ is a $\mathfrak{g}$-valued local
martingale and $\Gamma_{0}$ a $G$-valued $\mathcal{F}_{0}$-measurable random
variable (\cite[Proposition 1.9]{BSDE Lie groups}). This expresion provides
the bijection between $\mathfrak{g}$-valued local martingales and $\nabla
$-martingales on $G$ that we announced above.

\subsection{Geometric Brownian motion.}

Let $\left(  \mathbb{R}_{+},\cdot\right)  $ be the Abelian Lie group of
strictly positive real numbers endowed with the standard product. Its Lie
algebra is simply $\mathbb{R}$ and, for any $\xi\in\mathbb{R}$, the Lie
algebra exponential coincides with the standard exponential, that is $\exp
\xi=e^{\xi}$; consequently, the infinitesimal generator (right or
left-invariant) is
\[
\xi^{\mathbb{R}_{+}}(q)=\xi q,\text{ for any }q\in\mathbb{R}.
\]
Let $G=\mathbb{R}_{+}\times\overset{n}{\ldots}\times\mathbb{R}_{+}$ be the Lie
group constructed as the direct product of $n$ copies of $\left(
\mathbb{R}_{+},\cdot\right)  $. Its product map $\cdot:G\times G\rightarrow G$
is obviously $(a_{1},\ldots,a_{n})\cdot(b_{1},\ldots,b_{m})=(a_{1}b_{1}%
,\ldots,a_{n}b_{n})$, $a_{i}$, $b_{i}\in\mathbb{R}_{+}$ for any $i=1,\ldots
,n$, and its Lie algebra is $\mathfrak{g=}T_{1}\mathbb{R}_{+}\times\overset
{n}{\ldots}\times T_{1}\mathbb{R}_{+}\simeq\mathbb{R}\times\overset{n}{\ldots
}\times\mathbb{R}=\mathbb{R}^{n}$. Let $\{\xi_{i}=(0,\overset{i-1}{\ldots
},0,1,0,\ldots,0)\mid i=1,\ldots,n\}$ be the canonical basis of $\mathfrak{g}%
=\mathbb{R}^{n}$, $\mu=(\mu^{1},\ldots,\mu^{n})$, $\sigma=(\sigma^{1}%
,\ldots,\sigma^{n})\in\mathfrak{g}$ a couple of elements of $\mathfrak{g}$,
$B:\mathbb{R}_{+}\times\Omega\rightarrow\mathfrak{g}$ a $n$-dimensional
Brownian motion on some filtered probability space $\left(  \Omega
,P,\{\mathcal{F}_{t}\}_{t\in\mathbb{R}_{+}}\right)  $, and consider the
following Lie-Scheffers system on $G$
\begin{equation}
\delta\Gamma_{t}=\left(  \mu-\frac{1}{2}\sigma^{2}\right)  ^{G}(\Gamma
_{t})dt+\sum_{i=1}^{n}\sigma^{i}\xi_{i}^{G}(\Gamma_{t})\delta B_{t}^{i},
\label{eq 35}%
\end{equation}
where $\sigma^{2}=((\sigma^{1})^{2},\ldots,(\sigma^{n})^{2})$. Using
coordinates $\left(  q^{1},\ldots,q^{n}\right)  $ in $G$ we can rewrite
(\ref{eq 35}) as
\[
\delta q_{t}^{i}=\left(  \mu^{i}-\frac{1}{2}(\sigma^{i})^{2}\right)  q_{t}%
^{i}dt+\sigma^{i}q_{t}^{i}\delta B_{t}^{i},\quad i=1,\ldots,n,
\]
which may be rewritten in terms of It\^{o} integrals as
\begin{equation}
dq_{t}^{i}=\mu^{i}q_{t}^{i}dt+\sigma^{i}q_{t}^{i}dB_{t}^{i},~~i=1,\ldots,n.
\label{eq 36}%
\end{equation}
The solutions of the $n$-dimensional system of stochastic differential
equations (\ref{eq 36}) are usually referred to as the geometric Brownian
motion which is well-known for its use in the Black-Scholes theory of
derivatives pricing as a model for the time evolution of the prices of $n$
assets in a complete and arbitrage-free financial market.

The well-known solution of the differential equation~(\ref{eq 36}) can be
easily obtained by using the stochastic version of the Wei-Norman method that
we introduced in Section~\ref{wei norman method}. Indeed, let $q_{t}%
=\exp(a_{t}^{1}\xi_{1})\cdots\exp(a_{t}^{n}\xi_{n})$ be the solution of
(\ref{eq 36}) starting at $e=\left(  1,\ldots,1\right)  \in G$ as in the,
where $a^{i}:\mathbb{R}_{+}\times\Omega\rightarrow\mathbb{R}$ are real
semimartingales such that $a_{t=0}^{i}=0$ a.s. for any $i=1,\ldots,n$. Since
the Lie algebra $\mathfrak{g}$ of $G$ is Abelian, and~(\ref{eq 35}) is written
in Lie-Scheffers form
\[
\delta\Gamma_{t}=\sum_{i=1}^{l}\xi_{i}^{G}\left(  \Gamma_{t}\right)  \delta
X_{t}^{i}%
\]
by taking the noise semimartingale $X:=\left(  \left(  \mu^{1}-\frac
{(\sigma^{1})^{2}}{2}\right)  t+\sigma^{1}B_{t}^{1},\ldots,\left(  \mu
^{n}-\frac{(\sigma^{n})^{2}}{2}\right)  t+\sigma^{n}B_{t}^{n}\right)  $, the
equation (\ref{eq 33}) in the Wei-Norman method reduces to%
\[
\left(  \mu^{1}-(\sigma^{1})^{2}/2,\ldots,\mu^{n}-(\sigma^{n})^{2}/2\right)
t+\left(  \sigma^{1}B_{t}^{1},\ldots,\sigma^{n}B_{t}^{n}\right)  =\sum
_{i=1}^{n}\xi_{i}a_{t}^{i},
\]
which implies that $a_{t}^{i}=(\mu^{i}-(\sigma^{i})^{2}/2)t+\sigma^{i}%
B_{t}^{i}$ for any $i=1,\ldots,n$. Now, since the exponential map is given by
\[%
\begin{array}
[c]{rrl}%
\exp:\mathfrak{g} & \longrightarrow & G=\mathbb{R}_{+}^{n}\\
\xi=\sum_{i=1}^{n}\xi^{i}\xi_{i} & \longmapsto & \left(  \operatorname*{e}%
^{\xi^{1}},\ldots,\operatorname*{e}^{\xi^{n}}\right)
\end{array}
\]
where $\operatorname*{e}^{x}$ is the standard exponential function, we recover
the well-known result that the general solution $q_{t}$ of (\ref{eq 36})
starting at $q_{0}\in\mathbb{R}_{+}^{n}$ is%
\[
q_{t}=\left(  q_{0}^{1}\operatorname*{e}\nolimits^{(\mu^{1}-(\sigma^{1}%
)^{2}/2)t+\sigma^{1}B_{t}^{1}},\ldots,q_{0}^{n}\operatorname*{e}%
\nolimits^{(\mu^{n}-(\sigma^{n})^{2}/2)t+\sigma^{n}B_{t}^{n}}\right)  .
\]

\subsection{Brownian motion on reductive homogeneous spaces and symmetric
spaces.}

Let $G$ a Lie group and $H\subseteq G$ a closed subgroup. We say that the
homogeneous space $M=G/H$ is {\bfi{reductive}} if the Lie algebra
$\mathfrak{g}$ of $G$ may be decomposed into as a direct sum $\mathfrak{g}%
=\mathfrak{h}\oplus\mathfrak{m}$ where $\mathfrak{h}$ is the Lie algebra of
$H$ and $\mathfrak{m}$ is a subspace invariant under the action of
$\operatorname*{Ad}_{H}$. That is, $\operatorname*{Ad}\nolimits_{h}\left(
\mathfrak{m}\right)  \subseteq\mathfrak{m}$ for any $h\in H$ and,
consequently, $\left[  \mathfrak{h},\mathfrak{m}\right]  \subseteq
\mathfrak{m}$. Suppose now that the reductive homogeneous space $M$ is Riemann
manifold with Riemmanian metric $\eta$ and that the transitive action of $G$
leaves the metric $\eta$ invariant. We want to define Brownian motions on
$\left(  M,\eta\right)  $ by reducing a suitable process defined on $G$. The
notation and most of the results in this example, in addition to a
comprehensive exposition on homogeneous spaces, can be found in
\cite{helgason} and \cite{kobayashi-nomizu}. The reader is encouraged to check
with \cite{elworthy diffusion book} to learn more about the geometry of
homogeneous spaces in the stochastic context.

We start by recalling that a $M$-valued process $\Gamma$ is a Brownian motion
whenever%
\[
f(\Gamma)-f\left(  \Gamma_{0}\right)  -\frac{1}{2}\int\Delta(f)\left(
\Gamma_{s}\right)  ds
\]
is a real valued local semimartingale for any $f\in C^{\infty}(M)$, where
$\Delta$ denotes the Laplacian. As is widely known, the Laplacian is defined
as the trace of the Hessian associated to the Riemannian connection $\nabla$
of $\eta$. That is,%
\[
\Delta\left(  f\right)  \left(  m\right)  =\sum_{i=1}^{r}\left(
\mathcal{L}_{Y_{i}}\circ\mathcal{L}_{Y_{i}}-\nabla_{Y_{i}}Y_{i}\right)
(f)(m)
\]
where $\left\{  Y_{1},\ldots,Y_{r}\right\}  \subset\mathfrak{X}\left(
M\right)  $ is family or vector fields such that $\left\{  Y_{1}%
(m),\ldots,Y_{r}(m)\right\}  $ is an orthonormal basis of $T_{m}M$, $m\in M$.

Let $o\in M$ denote the equivalent class of $H$ in $M$. We have assumed that
$\left(  M,\eta\right)  $ is a Riemann manifold wiht a (left) $G$-invarinat
metric $\eta$. Since $\eta$ is $G$-invariant and $\Phi$ is transitive, the
only thing that really matters as far as the characterization of $\eta$ is
concerned is the symmetric bilinear form $\eta_{o}:T_{o}M\times T_{o}%
M\rightarrow T_{o}M$. It can be easily proved that there is a natural
one-to-one correspondence between the $G$-invariant Riemannian metrics $\eta$
on $M=G/H$ and the $\operatorname*{Ad}_{H}$-invariant positive definite
symmetric bilinear forms $B$ on $T_{o}M=\mathfrak{g}/\mathfrak{h}$
(\cite[Chapter X Proposition 3.1]{kobayashi-nomizu}). The correspondence is
given by%
\[
\eta\left(  \xi_{1}^{M},\xi_{2}^{M}\right)  =B\left(  T_{e}\pi\left(  \xi
_{1}\right)  ,T_{e}\pi\left(  \xi_{2}\right)  \right)  ,
\]
where $\xi_{1},\xi_{2}\in\mathfrak{g}$, $\pi:G\rightarrow G/H$ is the
canonical submersion, and $\xi^{M}\in\mathfrak{X}\left(  M\right)  $ denotes
the infinitesimal generator associated to $\xi\in\mathfrak{g}$. In addition,
if $M$ is reductive then the bilinear form $B$ may be regarded as defined on
$\mathfrak{m}$, $B:\mathfrak{m}\times\mathfrak{m}\rightarrow\mathbb{R}$, since
$T_{o}M$ is naturally isomorphic to $\mathfrak{m}$, which is an
$\operatorname*{Ad}_{H}$-invariant subspace of $\mathfrak{g}$. The Riemannian
connection $\nabla$ of the metric $\eta$ associated to such a bilinear form
$B$ is given by%
\begin{equation}
\nabla_{\xi_{1}^{M}}\xi_{2}^{M}=\frac{1}{2}\left[  \xi_{1}^{M},\xi_{2}%
^{M}\right]  +\left(  U\left(  \xi_{1},\xi_{2}\right)  \right)  ^{M},
\label{eq 45}%
\end{equation}
(\cite[Chapter X Theorem 3.3]{kobayashi-nomizu}). In this expression $\xi_{1}$
and $\xi_{2}$ belong to $\mathfrak{m}$ and $U:\mathfrak{m}\times
\mathfrak{m}\rightarrow\mathfrak{m}$ is the bilinear mapping defined by%
\[
2B\left(  U\left(  \xi_{1},\xi_{2}\right)  ,\xi_{3}\right)  =B\left(  \xi
_{1},\left[  \xi_{3,}\xi_{2}\right]  _{\mathfrak{m}}\right)  +B\left(  \left[
\xi_{3,}\xi_{1}\right]  _{\mathfrak{m}},\xi_{2}\right)  ,
\]
where $\left[  \cdot,\cdot\right]  _{\mathfrak{m}}$ is such that $\left[
\cdot,\cdot\right]  =\left[  \cdot,\cdot\right]  _{\mathfrak{h}}+\left[
\cdot,\cdot\right]  _{\mathfrak{m}}$ with $\left[  \cdot,\cdot\right]
_{\mathfrak{h}}\in\mathfrak{h}$ and $\left[  \cdot,\cdot\right]
_{\mathfrak{m}}\in\mathfrak{m}$. A consequence of (\ref{eq 45}) is that the
Laplacian $\Delta$ takes the expression $\Delta\left(  f\right)  \left(
m\right)  =\sum_{i=1}^{r}(\mathcal{L}_{\xi_{i}^{M}}\circ\mathcal{L}_{\xi
_{i}^{M}}+U\left(  \xi_{i},\xi_{i}\right)  ^{M})\left(  f\right)  (m)$, $m\in
M=G/K$, where $\{\xi_{1}^{M},\ldots,\xi_{r}^{M}\}$ is an orthonormal basis of
$T_{m}M$.

The most important examples of reductive homogeneous spaces are symmetric
spaces. In that case, $G$ is the connected component of the isometric group
$I(M)\subseteq\operatorname{Diff}(M)$ of the symmetric space $(M,\eta)$
containing $e=\operatorname*{Id}$. In order to identify the symmetric space
$\left(  M,\eta\right)  $ with a reductive space, take $o\in M$ a fixed point
and let $s$ a geodesic symmetry at $o$. Then the Lie group $G$ acts on $M$
transitively and, if $H$ denotes the isotropy group of $o$, $M$ is
diffeomorphic to $G/H$ (\cite[Chapter IV Theorem 3.3]{helgason}). Suppose that
$\dim\left(  G\right)  <\infty$ and let $\sigma:G\rightarrow G$ be the
involutive automorphism of $G$ defined by $\sigma\left(  g\right)  =s\circ
\Phi_{g}\circ s$ for any $g\in G$, where $\Phi:G\times M\rightarrow G$ denotes
as usual the left action of $G$ on $M$. It is a matter of fact that
$T_{e}\sigma:\mathfrak{g}\rightarrow\mathfrak{g}$ induces an involutive
automorphism of $\mathfrak{g}$. That is, $T_{e}\sigma\circ T_{e}%
\sigma=\operatorname*{Id}$ but $T_{e}\sigma\neq\operatorname*{Id}$. Let
$\mathfrak{h}$ and $\mathfrak{m}$ be the the eigenspaces of $\mathfrak{g}$
associated to the eigenvalues $1$ and $-1$ of $T_{e}\sigma$ respectively such
that $\mathfrak{g}=\mathfrak{h}\oplus\mathfrak{m}$. It can be checked that
$\mathfrak{h}$ is a Lie subalgebra of $\mathfrak{g}$,%
\[
\left[  \mathfrak{h},\mathfrak{h}\right]  \subseteq\mathfrak{h},\ \ \left[
\mathfrak{h},\mathfrak{m}\right]  \subseteq\mathfrak{m},\ \left[
\mathfrak{m},\mathfrak{m}\right]  \subseteq\mathfrak{h},
\]
and $\operatorname*{Ad}_{H}\left(  \mathfrak{m}\right)  \subseteq\mathfrak{m}$
(\cite[Chapter XI Proposition 2.1 and 2.2]{kobayashi-nomizu}). Morevoer, the
symmetric space $G/K$ has a unique affine connection $\nabla$ invariant under
the action of $G$. This is actually the Riemannian connection (\cite[Chapter
XI Theorem 3.3]{kobayashi-nomizu}) so that (\ref{eq 45}) reads%
\[
\nabla_{\xi_{1}^{M}}\xi_{2}^{M}=0
\]
for any pair of left-invariant vector fields $\xi_{1}^{M}$ and $\xi_{2}^{M} $.

Returning to the general case, let $\left\{  \xi_{1},\ldots,\xi_{r}\right\}  $
be a basis of $\mathfrak{m}$ such that $\left\{  T_{e}\pi\left(  \xi
_{1}\right)  \ldots,T_{e}\pi\left(  \xi_{r}\right)  \right\}  $ is an
orthonormal basis of $T_{o}(G/K)$ with respect to $\eta_{o}$ and let
$\{\xi_{1}^{G},\ldots,\xi_{r}^{G}\}\subset\mathfrak{X}\left(  G\right)  $ be
now the corresponding family of \textit{right-invariant} vector fields built
from $\left\{  \xi_{1},\ldots,\xi_{r}\right\}  $. Observe that $\{\xi_{1}%
^{M}(m),\ldots,\xi_{r}^{M}(m)\}$ is an orthonormal basis of $T_{m}(G/K)$ due
to the transitivity of the action and to the $G$-invariance of the metric
$\eta$. Consider now the Stratonovich stochastic differential equation%
\begin{equation}
\delta g_{t}=\sum_{i=1}^{r}\xi_{i}^{G}(g_{t})\delta B_{t}^{i}+\sum_{i=1}%
^{r}U\left(  \xi_{i},\xi_{i}\right)  ^{G}(g_{t})dt, \label{eq tan-nor 24}%
\end{equation}
where $\left(  B_{t}^{1},\ldots,B_{t}^{r}\right)  $ is a $\mathbb{R}^{r}%
$-valued Brownian motion. The stochastic system (\ref{eq tan-nor 24}) is by
definition $K$-invariant with respect to the natural \textit{right action}
$R:K\times G\rightarrow G$, $R_{k}\left(  g\right)  =gk$ for any $g\in G$ and
$k\in K$. In addition, it is straightforward to check that the projection
$\pi:G\rightarrow G/K$ send any right-invariant vector field $\xi^{G}%
\in\mathfrak{X}\left(  G\right)  $, $\xi\in\mathfrak{g}$, to the infinitesimal
generator $\xi^{M}\in\mathfrak{X}\left(  M\right)  $ of the $G$-action
$\Phi:G\times M\rightarrow M$. Hence (\ref{eq tan-nor 24}) projects to the
stochastic system
\begin{equation}
\delta\Gamma_{t}=\sum_{i=1}^{r}\xi_{i}^{M}(\Gamma_{t})\delta B_{t}^{i}%
+\sum_{i=1}^{r}U\left(  \xi_{i},\xi_{i}\right)  ^{M}\left(  \Gamma_{t}\right)
dt \label{eq tan-nor 25}%
\end{equation}
on $M$ by Proposition \ref{proposicion reduccion}. It is evident that the
solutions of (\ref{eq tan-nor 25}) have as a generator the second order
differential operator $\frac{1}{2}\sum_{i=1}^{r}(\mathcal{L}_{\xi_{i}^{M}%
}\circ\mathcal{L}_{\xi_{i}^{M}}+U\left(  \xi_{i},\xi_{i}\right)  ^{M})$ and
they are therefore Brownian motions.

\bigskip

\noindent\textbf{Acknowledgments} The authors are indebted to Jos\'{e} F.
Cari\~{n}ena and Eduardo Mart\'{\i}nez for their valuable comments and
suggestions. They acknowledge partial support from the French Agence National
de la Recherche, contract number JC05-41465. J.-A. L.-C. acknowledges support
from the Spanish Ministerio de Educaci\'{o}n y Ciencia grant number
BES-2004-4914. He also acknowledges partial support from MEC grant
BFM2006-10531 and Gobierno de Arag\'{o}n grant DGA-grupos consolidados
225-206. J.-P. O. has been partially supported by a \textquotedblleft Bonus
Qualit\'{e} Recherche" contract from the Universit\'{e} de Franche-Comt\'{e}.

\end{document}